\documentclass[oneside, 11pt]{amsart}
\usepackage[utf8]{inputenc}

\usepackage{amsfonts}
\usepackage{amsmath}
\usepackage{amssymb}
\usepackage{amsthm}
\usepackage[a4paper]{geometry}
\usepackage{fullpage}
\usepackage[dvipsnames]{xcolor}
\usepackage{mathtools}
\mathtoolsset{showonlyrefs}
\usepackage{caption}
\usepackage[colorlinks=true, linkcolor=blue, citecolor=blue]{hyperref}

\usepackage{enumerate}
\usepackage{tikz}
\usetikzlibrary{arrows.meta,calc,decorations.markings}
\tikzset{vtx/.style={circle, fill=black, inner sep=1.6pt}}

\usepackage{orcidlink}

\theoremstyle{plain}
\newtheorem{theorem}{Theorem}
\newtheorem{proposition}[theorem]{Proposition}
\newtheorem{lemma}[theorem]{Lemma}
\newtheorem{corollary}[theorem]{Corollary}

\newtheorem{remark}[theorem]{Remark}

\theoremstyle{definition}

\numberwithin{theorem}{section}
\numberwithin{equation}{section} 

\newcommand{\abs}[1]{\left\lvert #1 \right\rvert}
\newcommand{\one}{\mathbf{1}}

\newcommand{\Z}{\mathbb{Z}}

\newcommand{\N}{\mathbb{N}}

\newcommand{\calC}{\mathcal{C}}

\renewcommand{\P}{\mathbb P}
\newcommand{\eps}{\varepsilon}

\newcommand{\E}{\mathbb{E}}
\newcommand{\p}{\mathbb{P}}

\DeclareFontShape{OT1}{cmr}{m}{scit}{<->ssub*cmr/m/sc}{}

\makeatletter
\renewcommand\paragraph{\@startsection{paragraph}{4}{\z@}{.4\linespacing\@plus.4\linespacing}{-.5em}{\normalfont\bfseries}}
\makeatother

\title[Uniqueness in $k$-neighbor graphs]{Uniqueness and full coverage in $k$-neighbor graphs}
\author{David Coupier\orcidlink{0000-0002-1945-6337}}
\address{IMT Nord Europe, Institut Mines-T\'el\'ecom, Univ.\ Lille, F-59000 Lille, France}
\email{david.coupier@imt-nord-europe.fr}
\author{Beno\^it Henry}
\address{IMT Nord Europe, Institut Mines-T\'el\'ecom, Univ.\ Lille, F-59000 Lille, France}
\email{benoit.henry@imt-nord-europe.fr}
\author{Benedikt Jahnel\orcidlink{0000-0002-4212-0065}}
\address{Institut für Mathematische Stochastik, Technische Universit\"at Braunschweig \& Weierstrass Institute for Applied Analysis and Stochastics, Berlin}
\email{benedikt.jahnel@tu-braunschweig.de}
\author{Jonas K\"oppl\orcidlink{0000-0001-9188-1883}}
\address{Universität Münster, Germany.}
\email{jonas.koeppl@uni-muenster.de}

\date{\today}
\begin{document}

\begin{abstract}
We study the number of infinite clusters in the $2dp$-neighbor graph on $\Z^d$, $d\ge2$, in which every vertex independently selects $\lfloor 2dp\rfloor$ of its nearest neighbors uniformly at random, and one more with probability $2dp-\lfloor 2dp\rfloor$. The model has negatively correlated edges and, for $1\le2dp\le2d-1$, neither finite energy nor insertion or deletion tolerance. We prove that the undirected graph has at most one infinite cluster for every $p\in [0,1]$, and that in the planar case the same holds for forward and strong forward clusters of the directed graph when $p\in[1/2,1]$. We further exhibit a second phase transition: for $p\ge(d+1)/(2d)$ the open cluster of the origin is almost surely the whole lattice, and for $p>(d+1)/(2d)$ any two vertices of a box are connected inside the box of twice the size with probability exponentially close to one. In particular, the percolation function is continuous but not analytic on the supercritical interval $(p_c,1]$. Finally, we prove that the full-coverage phase transition at $p_f = (d+1)/2d$ is of order $2^d$. 
\end{abstract}

\maketitle

\section{Introduction}\label{sec:intro}
In the $k$-neighbor graph on $\Z^d$, every vertex independently selects $k$ of its $2d$ nearest neighbors uniformly at random and connects to them by directed edges~\cite{JKLT}. It is a lattice analogue of the $k$-nearest-neighbor graphs of continuum percolation~\cite{HaggMee96}. The continuous-parameter version introduced in~\cite{CHJK24}, the $2dp$-neighbor graph, lets every vertex select $\lfloor2dp\rfloor$ neighbors and, with probability $2dp-\lfloor 2dp\rfloor$, one more; thus every directed edge is open with probability $p$, and $p$ interpolates between the $k$-neighbor graphs. Forgetting orientations gives the undirected graph, and keeping only the edges selected by both endpoints gives the bidirectional graph. Percolation in these graphs, that is, the existence of an infinite cluster, has been studied in~\cite{JKLT,CHJK24}; in particular, the directed $2$-neighbor graph percolates in the plane~\cite{CHJK24}.

The constraint on the number of selected neighbors changes the nature of the model compared with classical Bernoulli percolation. The out-edges of a vertex are negatively correlated, and for $1\le2dp\le2d-1$ the model has neither finite energy nor insertion or deletion tolerance: since a vertex selects at least $\lfloor 2dp\rfloor$ and at most $\lceil 2dp\rceil$ neighbors, a single edge can in general be neither opened nor closed without changing other edges. For Bernoulli percolation, the proof of the uniqueness of the infinite cluster goes back to Aizenman, Kesten and Newman~\cite{aizenman1987uniqueness}. The  short argument for uniqueness of Burton and Keane~\cite{BurtonKeane89}, developed later, relies on finite energy, and insertion tolerance is also a key assumption for related results such as the indistinguishability of infinite clusters~\cite{LyonsSchramm99}. Uniqueness for models without insertion tolerance has recently been obtained by Panagiotis and Stauffer~\cite{PS26} for a class of monotone percolation models defined by local dynamics. In this paper we study the number of infinite clusters in the $2dp$-neighbor graphs and its consequences.

\paragraph{Main results.}
Our first result is uniqueness: for every $d\ge2$ and $p\in[0,1]$, the undirected $2dp$-neighbor graph has at most one infinite cluster (Theorem~\ref{thm_uni}). For $p>1/(2d)$ we prove this by a Burton--Keane argument in which the local modification rewires the out-sets of the vertices of a box, and for $p\le1/(2d)$ there is no infinite cluster. In the plane we also obtain uniqueness for the directed graph when $p\ge1/2$, both for forward clusters, in the sense that any two infinite forward clusters differ only in finitely many vertices, and for strong forward clusters, the sets of vertices that are mutually reachable. As for Bernoulli percolation~\cite{vdBK84,aizenman1987uniqueness}, uniqueness implies that the percolation function is continuous away from the critical point (Theorem~\ref{thm_con}).

Our second result is a second phase transition at $p=(d+1)/(2d)$, where the unique infinite cluster starts to cover the whole lattice. Below this value, finite clusters occur almost surely. At and above it,  the open cluster of the origin is all of $\Z^d$ almost surely. Strictly above the threshold this full coverage holds in a quantitative local form: with probability at least $1-c_1\exp(-c_2n)$, any two vertices of the box $\Lambda_n$ are connected inside $\Lambda_{2n}$ (Theorem~\ref{thm:full-coverage-local-uniqueness}). As a consequence, the percolation function of the undirected graph is continuous on the supercritical interval $(p_{\rm c},1]$ but equal to one on $[(d+1)/(2d),1]$ and smaller than one below, so it is not analytic there (Corollary~\ref{cor:non-analytic}). This contrasts with Bernoulli percolation on $\Z^d$, whose percolation function is analytic on the supercritical interval~\cite{GeorgakopoulosPanagiotis2023}. For Bernoulli percolation on transitive graphs, uniqueness of the infinite cluster is monotone in $p$~\cite{haggstrom1999monotonicity,schonmann1999stability}, which gives rise to a uniqueness transition; on $\Z^d$ it coincides with the percolation transition. The second transition of the $2dp$-neighbor graph is of a different nature: it concerns full coverage rather than uniqueness.

Our third result determines the order of this transition. The percolation function of the undirected graph is infinitely differentiable on $(1/2,1)$ except at $(d+1)/(2d)$. There, its derivatives up to order $2^d-1$ are continuous and vanish, whereas its derivative of order $2^d$ jumps from a negative value to zero. In particular, the probability that the origin's cluster is finite vanishes like $((d+1)/(2d)-p)^{2^d}$ as $p\uparrow(d+1)/(2d)$ (Theorem~\ref{thm:order}). Let us mention that non-analytic behavior of thermodynamic functions away from the critical point is known for example in the statistical mechanics of disordered systems as Griffiths singularities~\cite{griffiths1969nonanalytic,vanEnter2007griffiths}. In our model, both the non-analyticity and its finite order have a simple combinatorial origin: for $p\in[1/2,(d+1)/(2d))$, a finite cluster contains at least $2^d$ vertices that select exactly $d$ neighbors, one at each of its corners.

\paragraph{Outline.}
Section~\ref{sec:set} introduces the model and states our results, followed by a discussion of extensions and open problems. Section~\ref{sec:proofs} contains the proofs: preliminaries in Section~\ref{sec:prelim}, the Burton--Keane argument in Section~\ref{sec:BK}, full coverage and the sprinkling estimate in Section~\ref{sec:sprinkling}, the planar directed case in Section~\ref{sec:planar}, continuity in Section~\ref{sec:continuity} and the order of the full-coverage transition in Section~\ref{sec:order}.

\section{Setting and main results}\label{sec:set}
\subsection{Setting}
Throughout, $d\ge2$. We write $\mathcal N_o:=\{y\in\Z^d\colon\|y\|_{\ell^1}=1\}$ for the set of the $2d$ nearest neighbors of the origin~$o$, $\Lambda_n:=\{-n,\dots,n\}^d$ for the centered box of side length $2n+1$, and $\partial\Lambda_n:=\Lambda_n\setminus\Lambda_{n-1}$ for its inner boundary, $n\in\N$.

\paragraph{The directed $2dp$-neighbor graph.}
The model has a single parameter $p\in[0,1]$, from which we derive the two quantities
\begin{equation}
k:=\lfloor 2dp\rfloor\qquad\text{and}\qquad \eps:=2dp-k\in[0,1)\,.
\end{equation}
At each vertex $x\in\Z^d$, independently of all other vertices, we select $k$ of the $2d$ nearest neighbors of $x$ uniformly at random. In addition, the vertex performs an independent Bernoulli trial with success probability~$\eps$; if the trial is successful, one further neighbor is chosen uniformly among those not selected before. We call the set $A_x\subset x+\mathcal N_o$ of selected neighbors the \emph{out-set} of $x$, and we place a directed edge from $x$ to every vertex of $A_x$. Thus $|A_x|\in\{k,k+1\}$ with $|A_x|=k+1$ with probability $\eps$, and conditionally on $|A_x|=j$ the out-set $A_x$ is uniformly distributed among the $j$-element subsets of $x+\mathcal N_o$. In particular, every set of exactly $\lceil 2dp\rceil$ neighbors of $x$ is a possible out-set of $x$, and all of them are equally likely.

This procedure induces a probability measure $\P_p$ on the space $\Omega:=\{0,1\}^{\mathcal E}$ of configurations of directed edges, where
\begin{equation}
\mathcal E:=\{(x,y)\in(\Z^d)^2\colon\|x-y\|_{\ell^1}=1\}
\end{equation}
denotes the set of directed nearest-neighbor edges. A directed edge $e\in\mathcal E$ is \emph{open} (respectively \emph{closed}) in a configuration $\omega\in\Omega$ if $\omega(e)=1$ (respectively $\omega(e)=0$); thus $(x,y)$ is open if and only if $y\in A_x$.

The resulting model is the continuous-parameter extension, first introduced in~\cite{CHJK24}, of the $k$-neighbor graph studied in~\cite{JKLT}, and we refer to it as the \emph{directed $2dp$-neighbor graph}. Since $2dp=k+\eps$, the parameter $p$ and the pair $(k,\eps)$ determine each other, and we also write $\P_{(k,\eps)}:=\P_p$, which is sometimes more expressive. For instance, the model recovers the $k$-neighbor graph whenever $2dp=k$ is an integer, that is, $\P_{(k,0)}=\P_{k/(2d)}$ is the law of the directed $k$-neighbor graph; in dimension $d=2$, the parameters $p=1/4,1/2,3/4$ yield the directed $k$-neighbor graphs with $k=1,2,3$. By construction, the out-degree of every vertex is $k$ or $k+1$, and its expectation is
\begin{equation}
\E_p\big[|A_o|\big]=k(1-\eps)+(k+1)\eps=2dp\,,
\end{equation}
where $\E_p$ denotes the expectation with respect to $\P_p$. Consequently every directed edge is open with probability $p$, which motivates the parametrization. The out-edges of a fixed vertex are, however, not independent. They are negatively correlated in the sense that, for all $p\in(0,1)$ and all distinct $e,\hat e\in\mathcal E$ with the same initial vertex $x$,
\begin{equation}\label{eq:neg-corr}
\P_p(e\text{ is open}\mid\hat e\text{ is open})<p\,.
\end{equation}
Indeed, since all out-sets of a given size are equally likely, the probability that both $e$ and $\hat e$ are open equals $\E_p\big[|A_x|(|A_x|-1)\big]/(2d(2d-1))$. Using $\E_p\big[|A_x|^2\big]=(2dp)^2+\eps(1-\eps)$, one checks that~\eqref{eq:neg-corr} is equivalent to $\eps(1-\eps)<2dp(1-p)$. This is clear if $\eps=0$. Otherwise it follows by comparing $1/(\eps(1-\eps))=1/\eps+1/(1-\eps)$ with $1/(2dp(1-p))=1/(2dp)+1/(2d(1-p))$, since $\eps\le2dp$ and $1-\eps\le2d(1-p)$, and these two inequalities are not both equalities, as $2d>1$.

\begin{figure}[ht]
\input{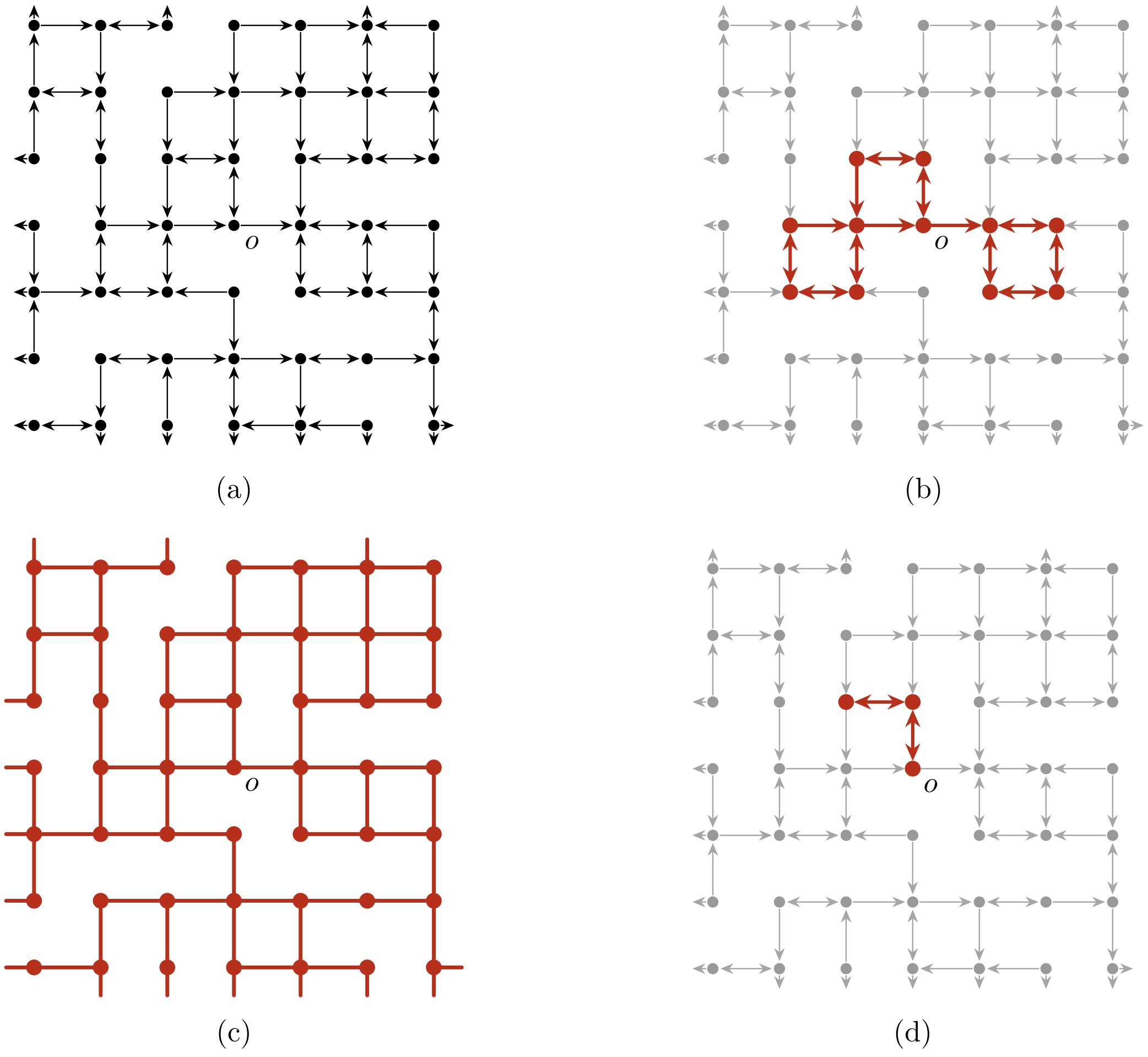}
\caption{A configuration of the directed $2$-neighbor graph on $\Z^2$, that is, $p=1/2$, in the window $\{-3,\dots,3\}^2$, and the clusters of the origin~$o$, which are highlighted in (b)--(d). (a) Every vertex of the window points to two of its four nearest neighbors, chosen uniformly at random. Double arrows indicate edges selected by both endpoints, and short arrows indicate edges leaving the window; edges selected by vertices outside the window are not drawn. (b) The forward cluster of~$o$, that is, the set of vertices that can be reached from $o$ along arrows, is finite here: no arrow leaves it. (c) In the undirected graph, every vertex of the window is connected to~$o$ inside the window, and the open cluster of~$o$ continues beyond the window. (d) The open cluster of~$o$ in the bidirectional graph consists of three vertices.}\label{fig_1}
\end{figure}

\paragraph{Undirected and bidirectional graphs.}
Two further models are derived from the directed one. In the \emph{undirected $2dp$-neighbor graph} an undirected edge $\{x,y\}$ is open if at least one of $(x,y)$ and $(y,x)$ is open. In the \emph{bidirectional $2dp$-neighbor graph} the edge $\{x,y\}$ is open if both $(x,y)$ and $(y,x)$ are open; see Figure~\ref{fig_1}.

\paragraph{Clusters and percolation.}
For $*\in\{\textsc{D},\textsc{U},\textsc{B}\}$ and $x,y\in\Z^d$ we write $x\rightsquigarrow^*y$ if there is a path of open edges from $x$ to $y$ in the directed (D), undirected (U) or bidirectional (B) $2dp$-neighbor graph, respectively; in the directed graph the path has to follow the orientation of its edges. Let $\calC^*(x):=\{y\in\Z^d\colon x\rightsquigarrow^* y\}$. Thus $\calC^\textsc{U}(x)$ and $\calC^\textsc{B}(x)$ are the open clusters of $x$, and $\calC^\textsc{D}(x)$ is the \emph{forward cluster} of $x$; see Figure~\ref{fig_1}(b)--(d). A stronger form of directed connectivity is captured by the \emph{strong forward cluster}
\begin{equation}
\calC^\textsc{S}(x):=\{y\in\Z^d\colon x\rightsquigarrow^\textsc{D}y\text{ and }y\rightsquigarrow^\textsc{D}x\}\subset\calC^\textsc{D}(x)\,.
\end{equation}
For $*\in\{\textsc{D},\textsc{S},\textsc{U},\textsc{B}\}$ we write $\{o\rightsquigarrow^*\infty\}:=\{|\calC^*(o)|=\infty\}$ and define the \emph{percolation function} and the \emph{critical threshold}
\begin{equation}
\theta^*_d(p):=\P_p(o\rightsquigarrow^*\infty)\qquad\text{and}\qquad p^*_{\rm c}(d):=\inf\{p\in[0,1]\colon\theta^*_d(p)>0\}\,.
\end{equation}
Each $\theta^*_d$ is nondecreasing by the monotone coupling of Lemma~\ref{lem:cont-coupling}, so there is at most one subcritical and one supercritical phase. A bidirectional open path can be traversed in both directions, so $\calC^\textsc{B}(x)\subset\calC^\textsc{S}(x)\subset\calC^\textsc{D}(x)\subset\calC^\textsc{U}(x)$, and therefore
\begin{equation}
p^\textsc{U}_{\rm c}(d)\le p^\textsc{D}_{\rm c}(d)\le p^\textsc{S}_{\rm c}(d)\le p^\textsc{B}_{\rm c}(d)\,.
\end{equation}
Moreover, $p^\textsc{U}_{\rm c}(d)\ge1/(2d)$, since the undirected $1$-neighbor graph does not percolate~\cite{JKLT}, and $p^\textsc{D}_{\rm c}(d)<1/2$ by~\cite[Corollary~2.4]{CHJK24}.

\paragraph{The number of infinite clusters.}
Our main focus lies on the number of infinite clusters. For $*\in\{\textsc{U},\textsc{B}\}$ we let $N^*$ denote the number of infinite open clusters of the respective graph. In the directed graph we have to be more careful, since forward clusters of different vertices may overlap without coinciding. For $*\in\{\textsc{D},\textsc{S}\}$ and vertices $x,y\in\Z^d$ with $|\calC^*(x)|=|\calC^*(y)|=\infty$ we write $x\sim^*y$ if $\calC^*(x)$ and $\calC^*(y)$ differ only in finitely many vertices, that is, if the symmetric difference $\calC^*(x)\triangle\calC^*(y)$ is finite. Since $A\triangle C\subset(A\triangle B)\cup(B\triangle C)$ for all sets $A,B,C$, this is an equivalence relation on the set of vertices with infinite cluster, and we let $N^*\in\{0,1,2,\dots\}\cup\{\infty\}$ be the number of its equivalence classes. We call $N^\textsc{D}$ the number of infinite forward clusters up to finite modifications. Thus $N^\textsc{D}\le1$ means precisely that any two infinite forward clusters differ only in finitely many vertices, that is, that there is at most one topological end. This is stronger than requiring that any two infinite forward clusters intersect: the forward cluster of $y$ may contain that of $x$ together with a further infinite directed path that $x$ does not reach; see Remark~\ref{rem:cont-ND}. Since two strong forward clusters are either disjoint or equal, and two disjoint infinite sets have infinite symmetric difference, $N^\textsc{S}$ is the number of infinite strong forward clusters. Each $N^*$, $*\in\{\textsc{D},\textsc{S},\textsc{U},\textsc{B}\}$, is invariant under translations and hence, by ergodicity, almost surely constant.

Topological ends are classically considered for random trees and forests, where a topological end is a class of infinite self-avoiding paths that eventually coincide, see for instance~\cite{HP98,GarciaSanchez26}. In the directed spanning forest on a planar Poisson point process, every point is linked to the nearest point with strictly larger abscissa. Thus every vertex has exactly one out-edge, forward clusters are infinite paths that are either disjoint or coalesce, and the analogue of $N^\textsc{D}$ is the maximal number of pairwise disjoint infinite paths; Coupier and Tran~\cite{CT13} showed by a Burton--Keane argument that it equals one. For supercritical oriented percolation on $\Z^d$, Grimmett and Hiemer~\cite{GH02} proved that, almost surely, any two vertices that are connected to infinity reach a common vertex by open directed paths, which is the analogue of the statement that any two infinite forward clusters intersect, and Wu and Zhang~\cite{wu2008geometrical} proved in two dimensions that any two infinite forward clusters differ only in finitely many vertices, which is the analogue of $N^\textsc{D}\le1$. Infinite strongly connected subgraphs, which correspond to infinite strong forward clusters, are studied for the continuum $k$-nearest-neighbor graph in~\cite{BalisterBollobas13}.

\subsection{Results}
Our first main result is the uniqueness of the infinite cluster.
\begin{theorem}[Uniqueness]\label{thm_uni}\hfill
\begin{enumerate}[(i)]
\item For all $d\ge2$ and $p\in[0,1]$ we have $N^\textsc{U}\in\{0,1\}$ almost surely.
\item For $d=2$, $*\in\{\textsc{D},\textsc{S}\}$ and $p\in[0,p^*_{\rm c}(2))\cup[1/2,1]$ we have $N^*\in\{0,1\}$ almost surely.
\end{enumerate}
\end{theorem}
Part (i) rests on a version of the Burton--Keane argument~\cite{BurtonKeane89}. Recall that for $1\le2dp\le2d-1$ the model has neither finite energy nor insertion or deletion tolerance, and the classical local modification, which opens all edges of a box, is in general not available. It is replaced by a rewiring of the out-sets of the vertices of a box that keeps the edges attaching the infinite clusters to the box (Section~\ref{sec:BK}). This argument works for every $p\in(1/(2d),1]$, on both sides of the full-coverage threshold $(d+1)/(2d)$ introduced below, while for $p\le1/(2d)$ there is no infinite cluster, by comparison with the undirected $1$-neighbor graph~\cite{JKLT}. Part (ii) relies on the planar structure exhibited in~\cite{CHJK24}: at $p=1/2$, every box is surrounded by an open directed cycle with high probability. For forward clusters, the bound $N^\textsc{D}\le1$ states that any two infinite forward clusters differ only in finitely many vertices (Lemma~\ref{lem:cont-planar}(i)), which is the form of uniqueness considered in~\cite{wu2008geometrical} for oriented percolation. Since $p^\textsc{D}_{\rm c}(2)<1/2$, Part (ii) gives $N^\textsc{D}=1$ for all $p\in[1/2,1]$. For $p<p^*_{\rm c}(2)$ trivially $N^*=0$; for $p\in[p^*_{\rm c}(2),1/2)$ the question remains open, see Section~\ref{sec:discussion}.

As for Bernoulli percolation~\cite{vdBK84,aizenman1987uniqueness}, uniqueness yields the continuity of the percolation function away from the critical point.
\begin{theorem}[Continuity]\label{thm_con}\hfill
\begin{enumerate}[(i)]
\item For $d\ge2$, the function $p\mapsto\theta^\textsc{U}_d(p)$ is continuous at every $p\in[0,1]\setminus\{p^\textsc{U}_{\rm c}(d)\}$.
\item For $d=2$ and $*\in\{\textsc{D},\textsc{S}\}$, the function $p\mapsto\theta^*_2(p)$ is continuous at every $p\in\big([0,p^*_{\rm c}(2))\cup[1/2,1]\big)\setminus\{p^*_{\rm c}(2)\}$.
\end{enumerate}
\end{theorem}
For forward and undirected clusters, right-continuity holds at every parameter, since percolation is a decreasing limit of events that depend on finitely many edges. Left-continuity at a supercritical parameter $p$ needs uniqueness at level $p$ in a form that connects the origin to the infinite clusters of all lower levels. In the directed case, this is exactly what the bound $N^\textsc{D}\le1$ provides, whereas the weaker statement that any two infinite forward clusters intersect is not enough (Remark~\ref{rem:cont-ND}). For strong forward clusters, right-continuity requires in addition the planar structure behind Theorem~\ref{thm_uni}(ii). This is why Part~(ii) is restricted to $p\ge1/2$. Whether $\theta^*_d$ is continuous at $p^*_{\rm c}(d)$ is open. For $*\in\{\textsc{D},\textsc{U}\}$, and for $*=\textsc{S}$ in the plane if $p^\textsc{S}_{\rm c}(2)\ge1/2$, this is equivalent to $\theta^*_d(p^*_{\rm c}(d))=0$, as in the Bernoulli case (Remark~\ref{rem:cont-critical}).

At $p=(d+1)/(2d)$ a second transition occurs. As soon as every vertex has at least $d+1$ out-edges, every forward cluster is infinite (Lemma~\ref{lem:forward-infinite}). Hence $\theta^\textsc{D}_d(p)=\theta^\textsc{U}_d(p)=1$ for $p\ge(d+1)/(2d)$, and in the undirected graph uniqueness takes the form of \emph{full coverage}: the open cluster of the origin is the whole lattice. Strictly above the threshold, the additional edges connect any two vertices of $\Lambda_n$ inside $\Lambda_{2n}$, with exponentially high probability.
\begin{theorem}[Full coverage and local uniqueness]\label{thm:full-coverage-local-uniqueness}
Let $d\ge2$.
\begin{enumerate}[(i)]
\item For every $p\ge(d+1)/(2d)$ we have $\P_p\big(\calC^\textsc{U}(o)=\Z^d\big)=1$.
\item For every $p>(d+1)/(2d)$ there are constants $c_1,c_2>0$, depending only on $d$ and $p$, such that for all $n\in\N$, 
\begin{equation}
	\P_p\big(\forall x,y \in \Lambda_n \colon x \rightsquigarrow^U y \text{ \rm via edges in } \Lambda_{2n} \big)
        \geq 
        1 - c_1 \exp(-c_2 n). 
\end{equation}
\end{enumerate}
\end{theorem}
In the strictly supercritical full-coverage regime it is thus exponentially unlikely that two vertices of a box fail to be connected locally. Part (i) follows from the uniqueness in Theorem~\ref{thm_uni}(i), since every vertex lies in an infinite open cluster, whereas Part (ii) is a consequence of a sprinkling estimate for locally dependent percolation models, Proposition~\ref{proposition:locally-dependent-sprinkling}, which may be of independent interest. It is inspired by the sprinkling results of~\cite{BHP2000,BT17} and by the Markov chain argument of Contreras and Hirschi~\cite{contrerasThesis}.

The full-coverage transition shows that the percolation function is not analytic in the supercritical regime, in contrast to Bernoulli percolation, for which it is analytic on $(p_{\rm c},1]$~\cite{GeorgakopoulosPanagiotis2023}; see Figure~\ref{fig_2}.
\begin{figure}[ht]
\centering
\begin{tikzpicture}[scale=1.15]
\def\pc{1.8}
\def\phalf{3}
\def\pend{4.5}
\draw[->,thick] (0,0) -- (6.3,0) node[right] {$p$};
\draw[->,thick] (0,-0.05) -- (0,4.5) node[left] {$\theta^{\textsc{U}}_d(p)$};
\draw (\pc,0.08)--(\pc,-0.08) node[below] {$p^{\textsc{U}}_{\rm c}(d)$};
\draw (\phalf,0.08)--(\phalf,-0.08) node[below] {$1/2$};
\draw (\pend,0.08)--(\pend,-0.08) node[below] {$(d+1)/(2d)$};
\draw (6,0.08)--(6,-0.08) node[below] {$1$};
\draw (-0.08,4)--(0.08,4) node[left] {$1\ $};
\draw[gray!50, dashed] (0,4)--(\pend,4);
\draw[very thick] (0,0)--(\pc,0);
\draw[very thick,domain=0:1,samples=200,variable=\u] plot ({\pc+(\pend-\pc)*\u},{4*sqrt(max(0,1-((1-\u)/(1-0.7*\u))^4))});
\draw[very thick] (\pend,4)--(6,4);
\fill (\pc,0) circle (1.5pt);
\fill (\pend,4) circle (1.5pt);
\node at (0.9,0.3) {\scriptsize $\theta^{\textsc{U}}_d(p)=0$};
\node at (5.25,3.7) {\scriptsize $\theta^{\textsc{U}}_d(p)=1$};
\end{tikzpicture}
\caption{Schematic shape of the percolation function $\theta^\textsc{U}_d$ of the undirected $2dp$-neighbor graph. It is continuous except possibly at $p^\textsc{U}_{\rm c}(d)$, infinitely differentiable on $(1/2,1)$ except at $(d+1)/(2d)$, equal to one on $[(d+1)/(2d),1]$ and smaller than one below $(d+1)/(2d)$, where $1-\theta^\textsc{U}_d(p)$ vanishes like $((d+1)/(2d)-p)^{2^d}$ (Theorem~\ref{thm_con}, Corollary~\ref{cor:non-analytic} and Theorem~\ref{thm:order}). The curve is drawn such that $1-\theta^\textsc{U}_d$ vanishes to order four at $(d+1)/(2d)$, as for $d=2$. The continuity at $p^\textsc{U}_{\rm c}(d)$ is only conjectured, and the shape of the curve below $1/2$ is illustrative.}\label{fig_2}
\end{figure}
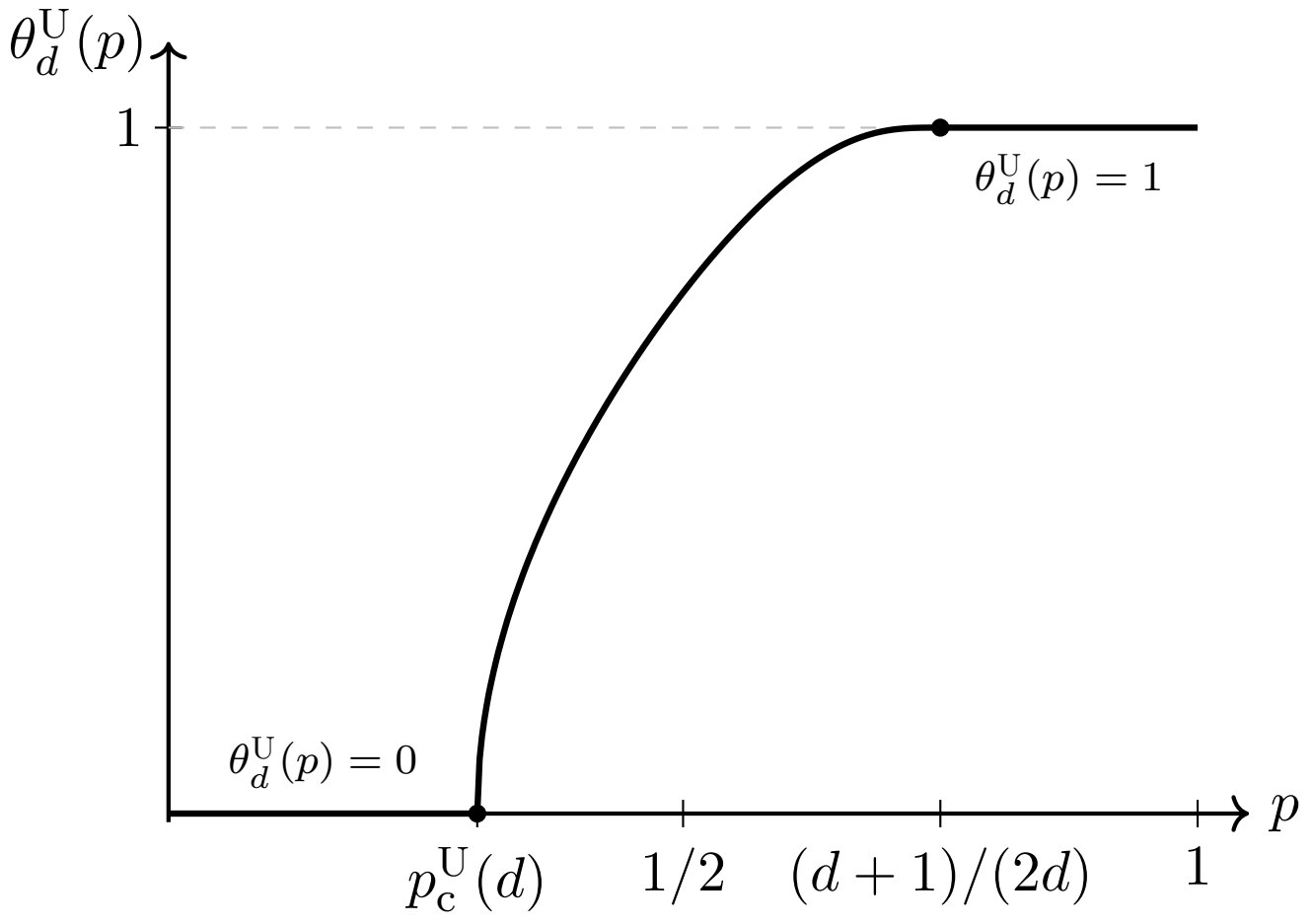
\begin{corollary}[Second phase transition]\label{cor:non-analytic}
Let $d\ge2$. Then, $\P_p\big(\calC^\textsc{U}(o)=\Z^d\big)=0$ and $\theta^\textsc{U}_d(p)<1$ for every $p<(d+1)/(2d)$, whereas $\P_p\big(\calC^\textsc{U}(o)=\Z^d\big)=\theta^\textsc{U}_d(p)=1$ for every $p\ge(d+1)/(2d)$. In particular, $\theta^\textsc{U}_d$ is continuous on $(p^\textsc{U}_{\rm c}(d),1]$, but it is not real-analytic on $(p^\textsc{U}_{\rm c}(d),1)$.
\end{corollary}
Corollary~\ref{cor:non-analytic} does not describe how $\theta^\textsc{U}_d$ approaches one. For $p\in[1/2,(d+1)/(2d))$, a vertex selects only $d$ neighbors with probability $d+1-2dp=2d((d+1)/(2d)-p)$, and $d+1$ neighbors otherwise. No out-edge leaves a finite open cluster. For each of the $2^d$ choices of signs, a vertex of the cluster that maximizes $\pm x_1\pm\dots\pm x_d$ over the vertices $x$ of the cluster has at most $d$ neighbors in the cluster, so it must select only $d$ neighbors. These $2^d$ \emph{corners} of the cluster are distinct (see the proof of Lemma~\ref{lem:corners}), so every finite cluster requires $2^d$ independent events, each of probability proportional to $(d+1)/(2d)-p$. Conversely, all vertices of the unit cube $\{0,1\}^d$ are corners, and the cube is the open cluster of the origin with probability of order $((d+1)/(2d)-p)^{2^d}$. One therefore expects $1-\theta^\textsc{U}_d(p)$ to vanish like $((d+1)/(2d)-p)^{2^d}$ as $p\uparrow(d+1)/(2d)$.
\begin{theorem}[Order of the full-coverage transition]\label{thm:order}
Let $d\ge2$.
\begin{enumerate}[(i)]
\item The function $\theta^\textsc{U}_d$ is infinitely differentiable on $(1/2,1)\setminus\{(d+1)/(2d)\}$ and $2^d-1$ times continuously differentiable on $(1/2,1)$.
\item There is a constant $c\in(0,\infty)$, depending only on $d$, such that
\begin{equation}
1-\theta^\textsc{U}_d(p)=c\,\Big(\frac{d+1}{2d}-p\Big)^{2^d}\Big(1+O\Big(\frac{d+1}{2d}-p\Big)\Big)\qquad\text{as }p\uparrow\frac{d+1}{2d}\,.
\end{equation}
\item The derivatives of $\theta^\textsc{U}_d$ of orders $1,\dots,2^d-1$ vanish at $(d+1)/(2d)$, and its derivative of order $2^d$ satisfies
\begin{equation}
\lim_{p\uparrow(d+1)/(2d)}\big(\theta^\textsc{U}_d\big)^{(2^d)}(p)=-(2^d)!\,c<0=\lim_{p\downarrow(d+1)/(2d)}\big(\theta^\textsc{U}_d\big)^{(2^d)}(p)\,.
\end{equation}
In particular, $\theta^\textsc{U}_d$ is $2^d-1$ times continuously differentiable on $(1/2,1)$ by~(i), but not $2^d$ times differentiable at $(d+1)/(2d)$.
\end{enumerate}
\end{theorem}
The derivative of order $2^d$ is thus the first derivative of $\theta^\textsc{U}_d$ that does not exist at $(d+1)/(2d)$: the full-coverage transition is of order four in the plane and of order eight in three dimensions. Turning the corner observation into Theorem~\ref{thm:order} requires control over large finite clusters, uniformly in $p$, since arbitrarily large clusters, such as boxes, have only $2^d$ corners. In Section~\ref{sec:order} we obtain this control by a Peierls argument, based on an elementary count of the possible boundaries of a finite cluster, and we deduce smoothness from Markov's inequality for derivatives of polynomials. Part~(i) also shows that $\theta^\textsc{U}_d$ is smooth on $(1/2,(d+1)/(2d))$; whether it is analytic there is open, see the following Section~\ref{sec:discussion}.

\subsection{Discussion}\label{sec:discussion}
Let us briefly mention some possible extensions of our results and related open problems for future investigations. This section can be skipped at first reading.
\paragraph{Other locally dependent models.}
Our arguments are not specific to the $2dp$-neighbor graph. They apply to locally dependent percolation models, in which every vertex independently selects a random set of neighbors according to a local law $\mu$ on the subsets of $\mathcal N_o$ (see Section~\ref{sec:sprinkling}), as long as $\mu$ is \emph{exchangeable}, that is, invariant under permutations of $\mathcal N_o$. Local laws that are supported on sets of at most one neighbor do not percolate, by comparison with the undirected $1$-neighbor graph~\cite{JKLT}. The Burton--Keane argument of Section~\ref{sec:BK} applies to every exchangeable local law that gives positive mass to some set of at least two neighbors: in Lemma~\ref{lem:rewiring}, the out-sets of the vertices of the box are replaced by sets of the largest size in the support of $\mu$, and this size is at least two. Together, these arguments give uniqueness of the infinite cluster for every exchangeable locally dependent percolation model on $\Z^d$, $d\ge2$, a class that contains many models without finite energy and without insertion or deletion tolerance. For exchangeable local laws that are supported on sets of at least $d+1$ neighbors, uniqueness takes the form of full coverage, as in Theorem~\ref{thm:full-coverage-local-uniqueness}(i), and if they also give positive mass to sets of at least $d+2$ neighbors, the domination of Lemma~\ref{lem:domination} and the sprinkling estimate of Proposition~\ref{proposition:locally-dependent-sprinkling} yield local uniqueness. Some parts of the proofs require less than exchangeability. In Section~\ref{sec:BK} it only serves to ensure that the rewired out-sets have positive probability, which is already the case if $\mu$ is \emph{quasi-exchangeable}, that is, if $\mu$ and its image under any permutation of $\mathcal N_o$ are mutually absolutely continuous.

\paragraph{Planar models.}
In the plane, the proof of Theorem~\ref{thm_uni}(ii) uses the model only through the crossing estimate~\eqref{eq:crossing}, the invariance of the model under translations and rotations, and monotonicity in $p$. We therefore expect Theorem~\ref{thm_uni}(ii) to hold also for the north-south/east-west (stick) and the corner models considered in~\cite{CHJK24}. For the $2dp$-neighbor graph itself, the restriction $p\ge1/2$ in Theorems~\ref{thm_uni}(ii) and~\ref{thm_con}(ii) comes from the crossing estimate, which is available at $p=1/2$ by~\cite{CHJK24} and hence, by monotonicity, for all $p\ge1/2$. If~\eqref{eq:crossing} held with $\P_{1/2}$ replaced by $\P_p$ for some $p\in(p^\textsc{D}_{\rm c}(2),1/2)$, the proofs of Sections~\ref{sec:planar} and~\ref{sec:continuity} would extend to all parameters in $[p,1]$. Whether $N^\textsc{D}\le1$ and $N^\textsc{S}\le1$ hold in dimensions $d\ge3$ is open.

\paragraph{The bidirectional graph.}
We do not know whether the bidirectional $2dp$-neighbor graph has at most one infinite cluster. If $N^\textsc{B}\le1$ almost surely, then $\theta^\textsc{B}_d$ is continuous on $[0,1]\setminus\{p^\textsc{B}_{\rm c}(d)\}$ by the argument for Theorem~\ref{thm_con}(i).

\paragraph{The critical point.}
As for Bernoulli percolation, the continuity of $\theta^\textsc{D}_d$ and $\theta^\textsc{U}_d$ at their critical points is equivalent to the absence of percolation at criticality (Remark~\ref{rem:cont-critical}), and we leave it open. For Bernoulli bond percolation on $\Z^d$, Kozma and Nitzan~\cite{KN24} reduced the absence of percolation at criticality to a conjectured inequality for percolation on finite graphs. Recently, a proof of this inequality was found by an AI system, and the resulting proof of the absence of percolation at criticality in all dimensions $d\ge2$ was verified formally in Lean~\cite{Leder26}; see~\cite{BlancRenaudie26,gouere2026proof} for a human-readable proof of the inequality and~\cite{BouRabeeLeder26} for an extension to site percolation on $\Z^d$, $d\ge3$. These arguments rely on the independence of the edges and on correlation inequalities, in particular on the Harris--FKG inequality. The latter fails for the $2dp$-neighbor graphs: for distinct $e,\hat e\in\mathcal E$ with the same initial vertex, the increasing events $\{e\text{ is open}\}$ and $\{\hat e\text{ is open}\}$ are negatively correlated by~\eqref{eq:neg-corr}, and since the out-sets of distinct vertices are independent, the same holds for two undirected edges with a common endpoint. Hence this route is not directly available for our model.

\paragraph{Smoothness of the percolation function.}
By Theorem~\ref{thm:order}, $\theta^\textsc{U}_d$ is infinitely differentiable on $(1/2,1)\setminus\{(d+1)/(2d)\}$. We do not know whether it is real-analytic on $(1/2,(d+1)/(2d))$. Our argument does not yield this, since the degrees of the polynomials in the proof of Proposition~\ref{prop:factorization} grow like $n^d$, whereas our bounds on them decay only exponentially in $n$. For $p<1/2$, vertices may select fewer than $d$ neighbors, so a finite cluster need not contain vertices that select exactly $d$ neighbors, and Lemma~\ref{lem:cut} is not available. We do not know whether $\theta^\textsc{U}_d$ is infinitely differentiable on $(p^\textsc{U}_{\rm c}(d),1/2]$. Note that the law of the out-sets is not differentiable in $p$ at the parameters $p=j/(2d)$, $j\in\N$, some of which lie in this interval, for instance $p=1/2$; at $p=1/2$, Proposition~\ref{prop:factorization} shows that $\theta^\textsc{U}_d$ has right derivatives of all orders.

\section{Proofs}\label{sec:proofs}
We first collect some general facts in Section~\ref{sec:prelim}. Section~\ref{sec:BK} contains the Burton--Keane argument for the undirected graph, which proves Theorem~\ref{thm_uni}(i), and Section~\ref{sec:sprinkling} the sprinkling argument behind Theorem~\ref{thm:full-coverage-local-uniqueness}. The planar directed case, Theorem~\ref{thm_uni}(ii), is treated in Section~\ref{sec:planar}, continuity, Theorem~\ref{thm_con}, in Section~\ref{sec:continuity}, and the order of the full-coverage transition, Theorem~\ref{thm:order}, in Section~\ref{sec:order}.

\subsection{Preliminaries}\label{sec:prelim}
Here we gather some basic properties and constructions that will be used later to prove our main results.
\paragraph{The monotone coupling.}
As in~\cite{CHJK24}, all measures $\P_p$, $p\in[0,1]$, can be realized on a common probability space, in such a way that the configurations are nondecreasing in $p$. Let $(\sigma_x,U_x)_{x\in\Z^d}$ be an independent family, where $\sigma_x$ is a uniformly distributed bijection from $\{1,\dots,2d\}$ onto $x+\mathcal N_o$ and $U_x$ is uniformly distributed on $[0,1)$ and independent of $\sigma_x$. We write $\P$ for the law of this family. Every directed edge $e=(x,y)\in\mathcal E$ receives the label
\begin{equation}\label{eq:cont-label}
T_e:=\frac{j-1+U_x}{2d}\in[0,1)\,,\qquad\text{where } j:=\sigma_x^{-1}(y)\,,
\end{equation}
and for $p\in[0,1]$ we define the configuration $\omega_p\in\Omega$ by $\omega_p(e):=\one\{T_e<p\}$. In words, $x$ opens its out-edges in the order $\sigma_x(1),\sigma_x(2),\dots$, the $j$-th one at level $(j-1+U_x)/(2d)$. Thus, with $k=\lfloor 2dp\rfloor$ and $\eps=2dp-k$, the out-neighbors of $x$ in $\omega_p$ are $\sigma_x(1),\dots,\sigma_x(k)$, together with $\sigma_x(k+1)$ if $U_x<\eps$, which is precisely the rule by which $\P_p$ selects the out-set of $x$. We order configurations coordinatewise, $\omega\le\omega'$ if $\omega(e)\le\omega'(e)$ for all $e\in\mathcal E$, and call an event $E\subset\Omega$ \emph{increasing} if $\omega\in E$ and $\omega\le\omega'$ imply $\omega'\in E$. The coupling satisfies the following basic properties, which we collect in Lemma~\ref{lem:cont-coupling} without proof for future reference.

\begin{lemma}\label{lem:cont-coupling}\hfill
\begin{enumerate}[(i)]
\item For every $p\in[0,1]$, the configuration $\omega_p$ has law $\P_p$.
\item If $p\le p'$, then $\omega_p\le\omega_{p'}$.
\item For every $e\in\mathcal E$, the label $T_e$ is uniformly distributed on $[0,1)$.
\end{enumerate}
\end{lemma}

For $*\in\{\textsc{D},\textsc{S},\textsc{U},\textsc{B}\}$, $x\in\Z^d$ and $p\in[0,1]$, $\calC^*_p(x)$ denotes the respective cluster of $x$ in the configuration $\omega_p$. By Lemma~\ref{lem:cont-coupling}, $\theta^*_d(p)=\P(|\calC^*_p(o)|=\infty)$ and $\calC^*_p(x)\subset\calC^*_{p'}(x)$ for $p\le p'$; in particular, $\theta^*_d$ is nondecreasing, so $\theta^*_d=0$ on $[0,p^*_{\rm c}(d))$ and $\theta^*_d>0$ on $(p^*_{\rm c}(d),1]$.

\paragraph{Infinite paths.}
Besides forward clusters we need the \emph{backward cluster}
\begin{equation}
\calC^{\leftarrow}(x):=\{y\in\Z^d\colon y\rightsquigarrow^\textsc{D}x\}\,,\qquad\text{so that}\qquad\calC^\textsc{S}(x)=\calC^\textsc{D}(x)\cap\calC^{\leftarrow}(x)\,,
\end{equation}
and $\calC^{\leftarrow}_p(x)$ denotes the backward cluster of $x$ in $\omega_p$. An \emph{infinite self-avoiding directed path starting at $x$} is a sequence $(y_j)_{j\ge0}$ of distinct vertices with $y_0=x$ and $(y_j,y_{j+1})\in\mathcal E$ for all $j$, and an \emph{infinite self-avoiding directed path ending at $x$} is a sequence $(y_j)_{j\ge0}$ of distinct vertices with $y_0=x$ and $(y_{j+1},y_j)\in\mathcal E$ for all $j$; such a path is open if all these edges are open.

\begin{lemma}\label{lem:cont-koenig}
Let $\omega\in\Omega$ and $x\in\Z^d$. If $|\calC^\textsc{D}(x)|=\infty$, then there is an infinite self-avoiding open directed path starting at $x$, and every vertex of this path has an infinite forward cluster. If $|\calC^{\leftarrow}(x)|=\infty$, then there is an infinite self-avoiding open directed path ending at $x$, and every vertex of this path has an infinite backward cluster.
\end{lemma}
\begin{proof}
Let $|\calC^\textsc{D}(x)|=\infty$. The finite self-avoiding open directed paths starting at $x$, partially ordered by extension, form a tree in which every element has at most $2d$ immediate successors. This tree is infinite, since every vertex of $\calC^\textsc{D}(x)$ is the endpoint of a shortest, hence self-avoiding, open directed path from $x$. By K\" onig's lemma the tree has an infinite branch, which is the required path, and every vertex of it reaches all later vertices. Since every vertex has at most $2d$ in-neighbors, the same argument applies to the finite self-avoiding open directed paths ending at $x$.
\end{proof}

\subsection{Uniqueness in the undirected graph}\label{sec:BK}
In this section we prove Theorem~\ref{thm_uni}(i). For $p\le1/(2d)$ there is no infinite cluster, and for every $p\in(1/(2d),1]$, on both sides of the full-coverage threshold $(d+1)/(2d)$, we use a Burton--Keane argument. The assumption $p>1/(2d)$, that is, $\lceil 2dp\rceil\ge2$, enters only in the rewiring of Lemma~\ref{lem:rewiring}. Throughout, \emph{open} refers to the undirected $2dp$-neighbor graph: the edge $\{x,y\}$ is open if $y\in A_x$ or $x\in A_y$.

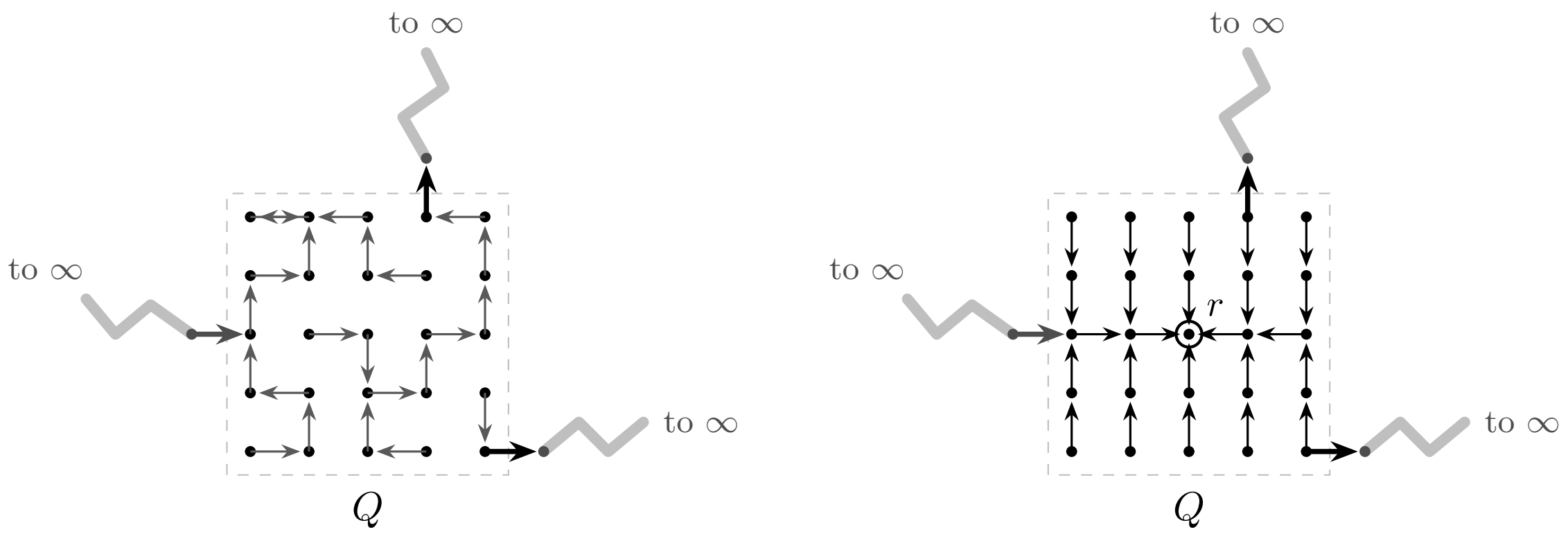
\begin{figure}[ht]
\centering
\begin{tikzpicture}[scale=0.56,
  v/.style={circle, fill=black, inner sep=1.05pt},
  w/.style={circle, fill=gray!60!black, inner sep=1.05pt},
  cl/.style={gray!50, line width=3pt, line cap=round, line join=round},
  keep/.style={-{Stealth[length=2.6mm]}, line width=1.4pt},
  tree/.style={-{Stealth[length=1.8mm]}, line width=0.6pt},
  old/.style={-{Stealth[length=1.8mm]}, line width=0.6pt, gray!70!black}]
  \newcommand{\clusters}{%
    \draw[cl] (-1,2) -- (-1.7,2.5) -- (-2.3,2.0) -- (-2.8,2.6);
    \draw[cl] (3,5) -- (2.6,5.7) -- (3.3,6.2) -- (3.0,6.8);
    \draw[cl] (5,0) -- (5.6,0.5) -- (6.1,0.0) -- (6.7,0.5);
    \node[gray!60!black, font=\footnotesize, anchor=south east] at (-2.6,2.7) {to $\infty$};
    \node[gray!60!black, font=\footnotesize, anchor=south] at (3.0,6.9) {to $\infty$};
    \node[gray!60!black, font=\footnotesize, anchor=west] at (6.8,0.5) {to $\infty$};
    \draw[gray!50, dashed] (-0.4,-0.4) rectangle (4.4,4.4);
    \node at (2,-1.0) {$Q$};
    \foreach \x in {0,...,4} \foreach \y in {0,...,4} \node[v] at (\x,\y) {};
    \node[w] at (-1,2) {}; \node[w] at (3,5) {}; \node[w] at (5,0) {};
    \draw[keep, gray!60!black] (-1,2) -- (-0.12,2);
    \draw[keep] (3,4) -- (3,4.88);
    \draw[keep] (4,0) -- (4.88,0);}
  \begin{scope}
    \clusters
    \draw[old] (0,0) -- (0.85,0);   \draw[old] (1,0) -- (1,0.85);
    \draw[old] (2,0) -- (2,0.85);   \draw[old] (3,0) -- (2.15,0);
    \draw[old] (0,1) -- (0,1.85);   \draw[old] (1,1) -- (0.15,1);
    \draw[old] (2,1) -- (2.85,1);   \draw[old] (3,1) -- (3,1.85);
    \draw[old] (4,1) -- (4,0.15);   \draw[old] (0,2) -- (0,2.85);
    \draw[old] (1,2) -- (1.85,2);   \draw[old] (2,2) -- (2,1.15);
    \draw[old] (3,2) -- (3.85,2);   \draw[old] (4,2) -- (4,2.85);
    \draw[old] (0,3) -- (0.85,3);   \draw[old] (1,3) -- (1,3.85);
    \draw[old] (2,3) -- (2,3.85);   \draw[old] (3,3) -- (2.15,3);
    \draw[old] (4,3) -- (4,3.85);   \draw[old] (0,4) -- (0.85,4);
    \draw[old] (1,4) -- (0.15,4);   \draw[old] (2,4) -- (1.15,4);
    \draw[old] (4,4) -- (3.15,4);
  \end{scope}
  \begin{scope}[xshift=14cm]
    \clusters
    \foreach \x in {0,...,4} {
      \draw[tree] (\x,0) -- (\x,0.85); \draw[tree] (\x,1) -- (\x,1.85);
      \draw[tree] (\x,4) -- (\x,3.15); \draw[tree] (\x,3) -- (\x,2.15);
    }
    \draw[tree] (0,2) -- (0.85,2); \draw[tree] (1,2) -- (1.85,2);
    \draw[tree] (4,2) -- (3.15,2); \draw[tree] (3,2) -- (2.15,2);
    \draw[thick] (2,2) circle (0.22);
    \node[font=\small] at (2.45,2.45) {$r$};
  \end{scope}
\end{tikzpicture}
\caption{Rewiring a box $Q$ (Lemma~\ref{lem:rewiring}). Left: three infinite clusters of the open graph outside $Q$ (thick gray) are attached to $Q$, one by an edge selected by the outside vertex and two by edges selected by vertices of $Q$ (bold arrows). Right: these attaching edges are kept, and every vertex of $Q$ other than $r$ also points to its neighbor on the way to $r$ in a spanning tree of $Q$ (thin arrows; further out-edges are not drawn). Now $Q$ is connected, and deleting it leaves at least three infinite clusters.}\label{fig:rewiring}
\end{figure}

The classical Burton--Keane argument~\cite{BurtonKeane89} consists of two parts. A deterministic counting argument shows that meeting points of three infinite clusters are rare in a large box, and a local modification shows that such meeting points would occur with positive probability if there were at least three infinite clusters. For Bernoulli percolation, or generally for models with finite energy, this modification opens all edges of a box. In our setting, this is in general not possible: a vertex selects at most $\lceil2dp\rceil$ neighbors, so out-edges cannot simply be added, and re-selecting the out-sets of the vertices of a box may cut infinite clusters off from the box. What remains is the freedom to rewire the fixed number of edges we get. The out-sets of the vertices of a box are independent of everything else, and all out-sets of a given size are equally likely. The main idea is to rewire the box to make it internally fully connected while also keeping the few out-edges that attach the infinite clusters to it, see Figure~\ref{fig:rewiring} below for a more detailed illustration.

\paragraph{Trifurcation boxes.}
A \emph{box} is a set $z+\Lambda_m$ with $z\in\Z^d$ and $m\ge0$. For a box $Q$, the \emph{open graph outside $Q$} is the graph on $\Z^d\setminus Q$ formed by the open edges with both endpoints outside $Q$. It depends only on the out-sets of the vertices outside $Q$, since an edge between two such vertices is open only if one of them selects the other. A cluster of the open graph outside $Q$ is \emph{attached} to $Q$ if it contains an endpoint of an open edge whose other endpoint lies in $Q$. We call $Q$ a \emph{trifurcation box} if
\begin{enumerate}[(a)]
\item any two vertices of $Q$ are joined by an open path inside $Q$, and
\item at least three infinite clusters of the open graph outside $Q$ are attached to $Q$.
\end{enumerate}
If $C$ is the open cluster containing the trifurcation box $Q$, then the clusters of the open graph outside $Q$ that are attached to $Q$ are precisely the connected components of $C\setminus Q$. Thus deleting a trifurcation box splits its cluster into at least three infinite pieces; we never need to close the edges between $Q$ and its surroundings.

\begin{figure}[ht]
\centering
\begin{tikzpicture}[scale=0.95,
  forest/.style={black, line width=0.8pt},
  tree/.style={black, line width=1.7pt},
  dead/.style={gray!60, line width=0.8pt},
  lat/.style={circle, fill=black, inner sep=0.95pt},
  leaf/.style={circle, fill=black, inner sep=1.5pt},
  deadend/.style={circle, fill=gray!60, inner sep=1.2pt},
  qv/.style={circle, fill=black, inner sep=2.4pt}]
  \def\panelframe{%
    \fill[gray!15] (-3,-3) rectangle (3,3);
    \fill[white] (-2.6,-2.6) rectangle (2.6,2.6);
    \draw[gray!70] (-3,-3) rectangle (3,3);
    \draw[gray!70, dashed] (-2.6,-2.6) rectangle (2.6,2.6);
    \node[gray!40!black, font=\footnotesize] at (-2.1,-2.3) {$\Lambda_n$};
  }
  \begin{scope}
    \panelframe
    \node at (0,3.4) {$F$};
    \foreach \cx/\cy/\lab/\lpos in {-1.1/0.5/Q_1/{(-1.62,-0.05)}, 1.0/-0.3/Q_2/{(1.0,0.38)}} {
      \fill[gray!35] (\cx-0.45,\cy-0.45) rectangle (\cx+0.45,\cy+0.45);
      \draw[forest] (\cx,\cy-0.3) -- (\cx,\cy+0.3);
      \foreach \j in {-0.3,0,0.3} \draw[forest] (\cx-0.3,\cy+\j) -- (\cx+0.3,\cy+\j);
      \foreach \i in {-0.3,0,0.3} \foreach \j in {-0.3,0,0.3} \node[lat] at (\cx+\i,\cy+\j) {};
      \node[font=\small] at \lpos {$\lab$};
    }
    \draw[forest] (-1.4,0.5) -- (-1.9,0.5) -- (-1.9,1.25) -- (-2.8,1.25);     
    \draw[forest] (-1.1,0.8) -- (-1.1,1.65) -- (-0.5,1.65) -- (-0.5,2.8);     
    \draw[forest] (-0.5,1.65) -- (0.35,1.65) -- (0.35,1.0);                    
    \draw[forest] (-0.8,0.5) -- (-0.15,0.5) -- (-0.15,-0.3) -- (0.7,-0.3);   
    \draw[forest] (-1.1,0.2) -- (-1.1,-0.75) -- (-1.75,-0.75) -- (-1.75,-1.45); 
    \draw[forest] (1.3,-0.3) -- (1.95,-0.3) -- (1.95,0.95) -- (2.8,0.95);     
    \draw[forest] (1.0,-0.6) -- (1.0,-1.5) -- (0.3,-1.5) -- (0.3,-2.8);       
    \foreach \p in {(-2.8,1.25),(-0.5,2.8),(2.8,0.95),(0.3,-2.8)} \node[leaf] at \p {};
    \foreach \p in {(0.35,1.0),(-1.75,-1.45)} \node[lat] at \p {};
  \end{scope}
  \draw[-{Stealth[length=2.2mm]}, gray!40!black] (3.35,0) -- (4.05,0);
  \begin{scope}[xshift=7.4cm]
    \panelframe
    \node at (0,3.4) {$\bar F$};
    \draw[dead] (-0.5,1.65) -- (0.35,1.65) -- (0.35,1.0);
    \draw[dead] (-1.1,0.5) -- (-1.1,-0.75) -- (-1.75,-0.75) -- (-1.75,-1.45);
    \foreach \p in {(0.35,1.0),(-1.75,-1.45)} \node[deadend] at \p {};
    \draw[tree] (-1.1,0.5) -- (-1.9,0.5) -- (-1.9,1.25) -- (-2.8,1.25);
    \draw[tree] (-1.1,0.5) -- (-1.1,1.65) -- (-0.5,1.65) -- (-0.5,2.8);
    \draw[tree] (-1.1,0.5) -- (-0.15,0.5) -- (-0.15,-0.3) -- (1.0,-0.3);
    \draw[tree] (1.0,-0.3) -- (1.95,-0.3) -- (1.95,0.95) -- (2.8,0.95);
    \draw[tree] (1.0,-0.3) -- (1.0,-1.5) -- (0.3,-1.5) -- (0.3,-2.8);
    \foreach \p in {(-2.8,1.25),(-0.5,2.8),(2.8,0.95),(0.3,-2.8)} \node[leaf] at \p {};
    \node[qv] at (-1.1,0.5) {};
    \node[qv] at (1.0,-0.3) {};
    \node[font=\small, anchor=north east] at (-1.17,0.43) {$q_1$};
    \node[font=\small, anchor=south west] at (1.07,-0.23) {$q_2$};
    \node[font=\small] at (1.45,-1.05) {$\mathcal T$};
  \end{scope}
\end{tikzpicture}
\caption{The construction in the proof of Lemma~\ref{lem:counting} for two trifurcation boxes $Q_1,Q_2\subset\Lambda_n$; the shaded frame is $\partial\Lambda_{n+1}$. The picture is schematic: paths are drawn as polylines, and only the component of $F$ that contains $Q_1$ and $Q_2$ is shown. Left: the spanning forest $F$ of $H$, which contains a spanning tree of each $Q_i$. Right: contracting these trees to single vertices $q_1$ and $q_2$ gives $\bar F$, whose component shown here is $K$. The smallest subtree $\mathcal T$ of $K$ containing the vertices of $K$ in $\partial\Lambda_{n+1}$ (bold) has all its leaves in $\partial\Lambda_{n+1}$, and $q_1$ and $q_2$ have degree three in it; branches that do not reach $\partial\Lambda_{n+1}$ (gray) do not belong to $\mathcal T$. Here $\mathcal T$ has four leaves and two vertices of degree three.}\label{fig:counting}
\end{figure}
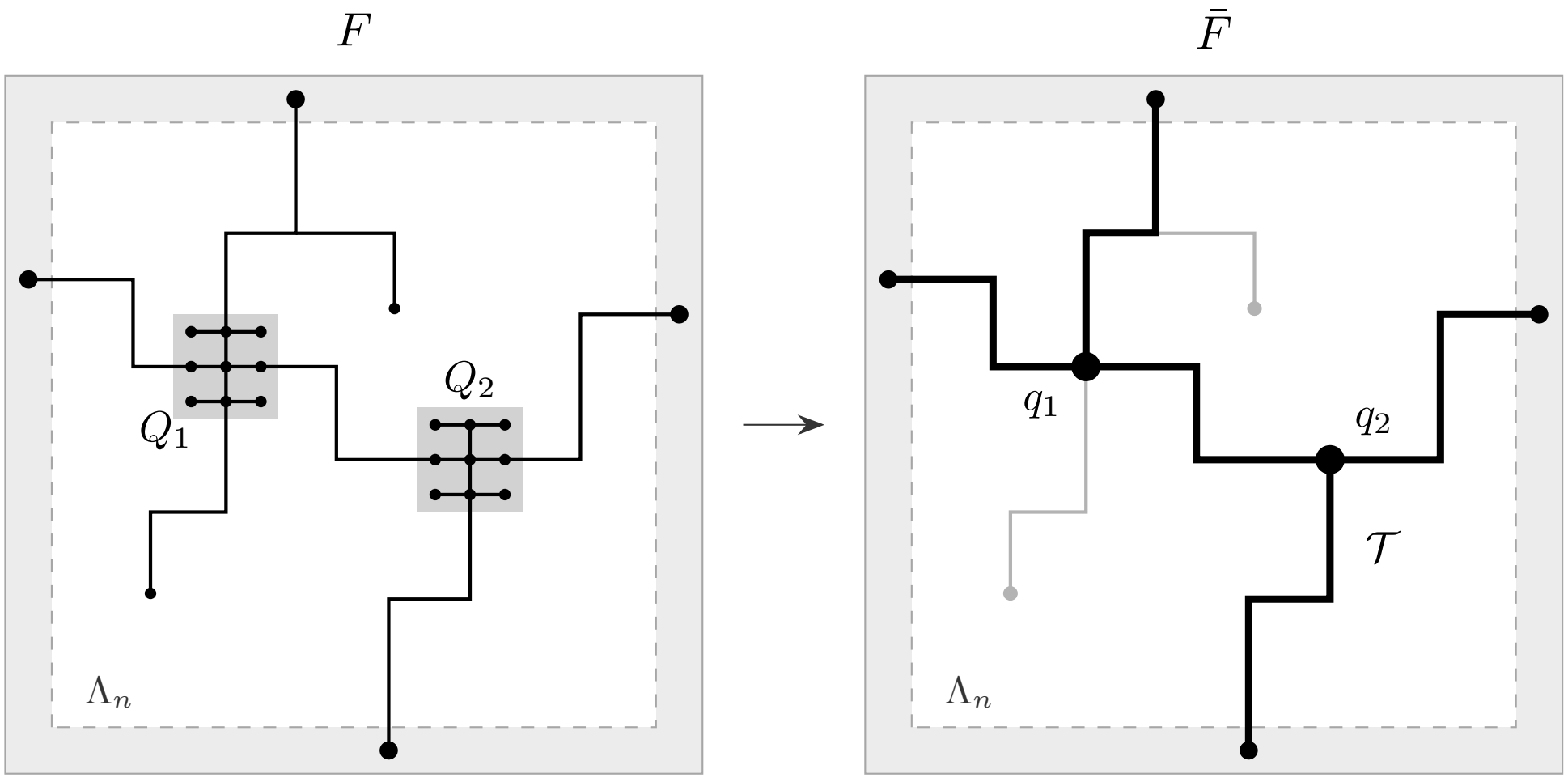

\begin{lemma}[Deterministic counting]\label{lem:counting}
Let $\omega\in\Omega$ and $n\in\N$. In every family of pairwise disjoint boxes contained in $\Lambda_n$, at most $|\partial\Lambda_{n+1}|$ boxes are trifurcation boxes.
\end{lemma}
\begin{proof}
This is the counting argument of Burton and Keane~\cite{BurtonKeane89}, see also~\cite[proof of Theorem~2.4]{HJ06}, with boxes in place of vertices. Let $Q_1,\dots,Q_r$ be the trifurcation boxes of the family. Let $H$ be the graph on $\Lambda_{n+1}$ whose edges are the open edges with at least one endpoint in $\Lambda_n$. By (a), and since $Q_i\subset\Lambda_n$, every $Q_i$ is connected in $H$. Since the boxes $Q_i$ are pairwise disjoint, the union of spanning trees of $Q_1,\dots,Q_r$ in $H$ is a forest, and we extend it to a spanning forest $F$ of $H$. Let $\bar F$ be the graph obtained from $F$ by contracting each spanning tree associated to $Q_i$ to a single vertex~$q_i$; see Figure~\ref{fig:counting}. Contracting the edges of pairwise disjoint subtrees of a forest creates no cycles, so $\bar F$ is again a forest.

For $i\in\{1,\dots,r\}$ let $\bar F-q_i$ denote the forest obtained by deleting the vertex $q_i$ from $\bar F$. We claim that the component of $\bar F$ containing $q_i$ contains vertices $b_1,b_2,b_3\in\partial\Lambda_{n+1}$ that lie in three distinct components of $\bar F-q_i$. In particular, at least three distinct components of $\bar F-q_i$ that are contained in the component of $\bar F$ containing $q_i$ contain a vertex of $\partial\Lambda_{n+1}$. To prove the claim, let $I_1,I_2,I_3$ be distinct infinite clusters of the open graph outside $Q_i$ attached to $Q_i$ by open edges $\{v_a,w_a\}$ with $v_a\in Q_i$ and $w_a\in I_a$. Since $I_a$ is infinite and connected, it contains a path from $w_a$ to a vertex $b_a\in\partial\Lambda_{n+1}$ whose other vertices lie in $\Lambda_n$; together with $\{v_a,w_a\}$ this is a path in $H$ from $Q_i$ to $b_a$. Hence $b_1,b_2,b_3$ lie in the component of $F$ containing $Q_i$, and thus in the component of $\bar F$ containing $q_i$; they are distinct, since $I_1,I_2,I_3$ are disjoint, and they are not contracted, since all boxes lie in $\Lambda_n$. Every path in $H$ from $b_a$ to $b_{a'}$ with $a\ne a'$ passes through $Q_i$, since its edges are open and $I_a$ and $I_{a'}$ are different clusters of the open graph outside $Q_i$. Consequently, the path from $b_a$ to $b_{a'}$ in the tree $\bar F$ passes through $q_i$: otherwise, replacing every vertex $q_l$, $l\ne i$, on it by a path in the spanning tree of $Q_l$ would give a path in $H$ from $b_a$ to $b_{a'}$ that avoids $Q_i$. Hence $b_1,b_2,b_3$ lie in three distinct components of $\bar F-q_i$, which proves the claim.

Now let $K$ be a component of $\bar F$ that contains some of the vertices $q_i$, and let $\mathcal T$ be the smallest subtree of $K$ that contains all vertices of $K$ in $\partial\Lambda_{n+1}$ (bold in Figure~\ref{fig:counting}); by the claim, there are at least three such vertices. By minimality of $\mathcal T$, all leaves of $\mathcal T$ lie in $\partial\Lambda_{n+1}$: if a leaf of $\mathcal T$ did not lie in $\partial\Lambda_{n+1}$, then deleting it together with its incident edge would give a smaller subtree of $K$ that still contains all vertices of $K$ in $\partial\Lambda_{n+1}$. Let $q_i$ be a vertex of $K$, and let $b_1,b_2,b_3$ be as in the claim; they lie in $K$. For $a\ne a'$, the path in $\bar F$ from $b_a$ to $b_{a'}$ lies in $\mathcal T$, since $\mathcal T$ is a subtree of $K$ containing $b_a$ and $b_{a'}$, and it passes through $q_i$, since $b_a$ and $b_{a'}$ lie in distinct components of $\bar F-q_i$. Hence $q_i$ belongs to $\mathcal T$. Moreover, the paths in $\mathcal T$ from $q_i$ to $b_1$, $b_2$ and $b_3$ begin with three distinct edges, since they lead into the three distinct components of $\bar F-q_i$ that contain $b_1$, $b_2$ and $b_3$; thus $q_i$ has degree at least three in $\mathcal T$. For a vertex $v$ of $\mathcal T$, which we denote by $v\in\mathcal T$, let $\deg_{\mathcal T}(v)$ be its degree in $\mathcal T$. Since $\mathcal T$ is a finite tree with at least two vertices, it has one edge fewer than vertices, so that $\sum_{v\in\mathcal T}(\deg_{\mathcal T}(v)-2)=-2$. In this sum every leaf contributes $-1$ and every vertex of degree two contributes $0$, and no vertex has degree zero. Hence
\begin{equation}
\abs{\{v\in\mathcal T\colon v\text{ is a leaf}\}}=2+\sum_{v\in\mathcal T\colon\deg_{\mathcal T}(v)\ge3}\big(\deg_{\mathcal T}(v)-2\big)\ge\abs{\{v\in\mathcal T\colon\deg_{\mathcal T}(v)\ge3\}}+2\,.
\end{equation}
Since all leaves of $\mathcal T$ lie in $\partial\Lambda_{n+1}$ and every vertex $q_i$ of $K$ has degree at least three in $\mathcal T$, it follows that $K$ contains fewer vertices $q_i$ than vertices in $\partial\Lambda_{n+1}$. Summing over all components $K$ of $\bar F$ that contain some of the vertices $q_i$, which are disjoint, we obtain $r\le|\partial\Lambda_{n+1}|$.
\end{proof}

The local modification is the following rewiring of a box, illustrated in Figure~\ref{fig:rewiring}.
\begin{lemma}[Rewiring a box]\label{lem:rewiring}
Let $p\in(1/(2d),1]$, let $Q$ be a box and let $j\in\{2,3\}$. Let $E$ be the event that at least $j$ infinite clusters of the open graph outside $Q$ are attached to $Q$, and let $F$ be the event that any two vertices of $Q$ are joined by an open path inside $Q$. There is a constant $c>0$, depending only on $d$, $p$ and $|Q|$, such that for every event $G$ that depends only on the out-sets of the vertices outside $Q$,
\begin{equation}\label{eq:rewiring}
\P_p(F\cap E\cap G)\ge c\,\P_p(E\cap G)\,.
\end{equation}
\end{lemma}
 In words, if $j$ infinite clusters are attached to $Q$, then $Q$ can be connected internally by rewiring the out-sets of its vertices, at the cost of the factor~$c$. Note that for $j=3$, the event $F\cap E$ is precisely the event that $Q$ is a trifurcation box.

 The main idea behind the rewiring in $Q$ is as follows. We can choose a spanning tree for $Q$ and then root it at some vertex $r \in Q$. Now each $v \in Q \setminus \{r\}$ only has to open the edge to its unique parent to make $Q$ internally connected while still satisfying the constraint. To make sure that the infinite clusters outside of $Q$ are then still connected to $Q$ we have to be a bit careful with choosing the root $r$, but the mechanism remains the same. 
 
\begin{proof}
Write $\xi:=(A_y)_{y\notin Q}$ and $\eta:=(A_v)_{v\in Q}$. Under $\P_p$ they are independent, and $\eta$ takes finitely many values, in all of which every out-set has at most $\lceil2dp\rceil$ elements.

We first show that for every configuration in $E$ in which the vertices of $Q$ have out-sets of at most $\lceil2dp\rceil$ elements there are sets $A'_v\subset v+\mathcal N_o$, $v\in Q$, each with exactly $\lceil2dp\rceil$ elements, such that replacing $A_v$ by $A'_v$ for all $v\in Q$ produces a configuration in $F\cap E$. All choices below are made by a fixed rule, for instance the first admissible one in a fixed enumeration. Let $I_1,\dots,I_j$ be distinct infinite clusters of the open graph outside $Q$, attached to $Q$ by open edges $\{v_a,w_a\}$ with $v_a\in Q$ and $w_a\in I_a$. If $v_a\in A_{w_a}$, this edge stays open whatever the vertices of $Q$ select. Otherwise $w_a\in A_{v_a}$, and we record $w_a$ in a set $D_{v_a}$. In this way we obtain sets $D_v\subset A_v\setminus Q$ of out-edges that we keep, with $\sum_{v\in Q}|D_v|\le j\le3$. Since $|D_v|\le|A_v|\le\lceil2dp\rceil$ and $\lceil 2dp\rceil\ge2$ by the assumption $p>1/(2d)$, at most one vertex $v\in Q$ satisfies $|D_v|=\lceil 2dp\rceil$. Let $r$ be this vertex if it exists, and an arbitrary vertex of $Q$ otherwise. Fix a spanning tree of $Q$, viewed as a subgraph of the lattice. We will root this tree at $r$ and then use the generation structure of this tree to make sure that no vertex's choices exceed its allowed budget. For $v\in Q\setminus\{r\}$ let $v'$ be the parent of $v$, i.e., the neighbor of $v$ on the path from $v$ to $r$ in this tree. Choose $A'_r\supset D_r$ and $A'_v\supset D_v\cup\{v'\}$ for $v\ne r$, all of size $\lceil 2dp\rceil$. This is possible, since $v'\in Q$ while $D_v\cap Q=\emptyset$, and $|D_v|+1\le\lceil2dp\rceil\le2d$ for $v\ne r$. In the new configuration the tree edges $\{v,v'\}$ are open, so any two vertices of $Q$ are joined by an open path inside $Q$; the edges $\{v_a,w_a\}$ are still open; and the open graph outside $Q$ is unchanged. Hence the new configuration lies in $F\cap E$.

Now fix an enumeration of the possible values of $\eta$. On the event that some possible value $\eta^*$ satisfies $(\xi,\eta^*)\in E$, choose the first such $\eta^*$ and let $\phi(\xi)$ be the family $(A'_v)_{v\in Q}$ constructed above from $(\xi,\eta^*)$; as $\phi$ takes finitely many values, it is a measurable function of $\xi$. Every family of out-sets of size $\lceil2dp\rceil$ of the vertices of $Q$ has the same probability $c>0$, namely $c=\big(\eps\binom{2d}{k+1}^{-1}\big)^{|Q|}$ if $\eps>0$ and $c=\binom{2d}{k}^{-|Q|}$ if $\eps=0$. If $\eta=\phi(\xi)$, then the configuration lies in $F\cap E$. Since $G$ depends only on $\xi$, and $\eta$ is independent of $\xi$,
\begin{align*}
\P_p(F\cap E\cap G)&\ge\P_p\big(G\cap\{\phi(\xi)\text{ is defined and }\eta=\phi(\xi)\}\big)\\
&=c\,\P_p\big(G\cap\{\phi(\xi)\text{ is defined}\}\big)\ge c\,\P_p(E\cap G)\,.\qedhere
\end{align*}
\end{proof}

\begin{proposition}\label{prop:BK}
Let $d\ge2$ and $p\in(1/(2d),1]$. Then $N^\textsc{U}\le1$ $\P_p$-almost surely.
\end{proposition}
\begin{proof}
\emph{Step 1: Trifurcation boxes do not occur.} Fix $m\ge0$. The boxes $(2m+1)z+\Lambda_m$, $z\in\Z^d$, are pairwise disjoint, and those with $\|z\|_\infty\le n$ are contained in $\Lambda_{(2m+1)n+m}$. Taking expectations in Lemma~\ref{lem:counting} and using translation invariance, we obtain for all $n\in\N$
\begin{equation}
(2n+1)^d\,\P_p(\Lambda_m\text{ is a trifurcation box})\le\big|\partial\Lambda_{(2m+1)n+m+1}\big|\le2d\,\big(2(2m+1)n+2m+3\big)^{d-1}\,.
\end{equation}
Letting $n\to\infty$ shows that $\Lambda_m$ is a trifurcation box with probability zero. Now let $E$ and $F$ be the events of Lemma~\ref{lem:rewiring} for $Q=\Lambda_m$ and $j=3$, so that $F\cap E$ is the event that $\Lambda_m$ is a trifurcation box. We apply~\eqref{eq:rewiring} with $G=\Omega$, the sure event, which trivially depends only on the out-sets of the vertices outside $\Lambda_m$, and obtain
\begin{equation}
c\,\P_p(E)\le\P_p(F\cap E)=\P_p(\Lambda_m\text{ is a trifurcation box})=0\,.
\end{equation}
Hence, almost surely, at most two infinite clusters of the open graph outside $\Lambda_m$ are attached to $\Lambda_m$.

\emph{Step 2: Infinite clusters meeting a box.} Let $C$ be an infinite open cluster that meets $\Lambda_m$. Since $C$ is connected, every connected component of $C\setminus\Lambda_m$ is a cluster of the open graph outside $\Lambda_m$ attached to $\Lambda_m$. There are finitely many such components, since each contains a neighbor of $\Lambda_m$, so at least one of them is infinite. Consequently, distinct infinite open clusters meeting $\Lambda_m$ give rise to distinct infinite clusters of the open graph outside $\Lambda_m$ attached to $\Lambda_m$. Moreover, the infinite clusters of the open graph outside $\Lambda_m$ that are not attached to $\Lambda_m$ are precisely the infinite open clusters disjoint from $\Lambda_m$.

\emph{Step 3: Conclusion.} By ergodicity, $N^\textsc{U}$ is almost surely equal to a constant $N\in\{0,1,2,\dots\}\cup\{\infty\}$. If $N\ge3$, then for $m$ large enough, with positive probability at least three infinite open clusters meet $\Lambda_m$; by Step~2, at least three infinite clusters of the open graph outside $\Lambda_m$ are then attached to $\Lambda_m$, which contradicts Step~1. Suppose that $N=2$, and choose $m$ such that with positive probability both infinite open clusters meet $\Lambda_m$; call this event $B$. By Steps 1 and 2, almost surely on $B$ the open graph outside $\Lambda_m$ has exactly two infinite clusters, and both are attached to $\Lambda_m$. Let $G$ be the event that the open graph outside $\Lambda_m$ has at most two infinite clusters, which depends only on the out-sets of the vertices outside $\Lambda_m$, and let $E$ and $F$ now be the events of Lemma~\ref{lem:rewiring} for $Q=\Lambda_m$ and $j=2$. Almost surely on $B$, both $E$ and $G$ occur, so~\eqref{eq:rewiring} gives
\begin{equation}
\P_p(F\cap E\cap G)\ge c\,\P_p(E\cap G)\ge c\,\P_p(B)>0\,.
\end{equation}
On $F\cap E\cap G$, the box $\Lambda_m$ is connected by open paths inside it, and the open graph outside $\Lambda_m$ has exactly two infinite clusters, both attached to $\Lambda_m$. By Step 2, there is then no infinite open cluster disjoint from $\Lambda_m$, and every infinite open cluster meets, and hence contains, $\Lambda_m$. Thus $N^\textsc{U}=1$ on $F\cap E\cap G$, which contradicts $N=2$. Hence $N\le1$.
\end{proof}

\begin{proof}[Proof of Theorem~\ref{thm_uni}(i)]
For $p\le1/(2d)$ we have $\omega_p\le\omega_{1/(2d)}$ by Lemma~\ref{lem:cont-coupling}, and the undirected $1$-neighbor graph, that is, the undirected $2dp$-neighbor graph with $p=1/(2d)$, does not percolate~\cite{JKLT}; hence $N^\textsc{U}=0$ almost surely. For $p\in(1/(2d),1]$ the claim is Proposition~\ref{prop:BK}.
\end{proof}

\subsection{Full coverage and local uniqueness}\label{sec:sprinkling}
In this section we prove Theorem~\ref{thm:full-coverage-local-uniqueness} and Corollary~\ref{cor:non-analytic}. The starting point is a deterministic observation.

\begin{lemma}\label{lem:forward-infinite}
Let $\omega\in\Omega$ be a configuration in which every vertex has at least $d+1$ open out-edges. Then every vertex has an infinite forward cluster. In particular, for $p\ge(d+1)/(2d)$, $\P_p$-almost surely every vertex has an infinite forward cluster, and $\theta^\textsc{D}_d(p)=\theta^\textsc{U}_d(p)=1$.
\end{lemma}
\begin{proof}
Suppose that $\calC^\textsc{D}(x)$ is finite for some $x$, and let $v$ be a vertex of $\calC^\textsc{D}(x)$ for which the sum $v_1+\dots+v_d$ of its coordinates is maximal. All out-neighbors of $v$ belong to $\calC^\textsc{D}(x)$, so none of them is of the form $v+e_i$, where $e_1,\dots,e_d$ denote the unit vectors. Hence all out-neighbors of $v$ are among $v-e_1,\dots,v-e_d$, contradicting the assumption that $v$ has at least $d+1$ of them. For $p\ge(d+1)/(2d)$ we have $k\ge d+1$, so $\P_p$-almost surely every vertex has at least $d+1$ open out-edges.
\end{proof}
A similar argument is used in~\cite{JKLT} for the directed $k$-neighbor graph with $k\ge d+1$. Above the threshold $(d+1)/(2d)$, some vertices have additional out-edges, which act as a sprinkling on top of an everywhere percolating graph. This sprinkling is neither a family of independent edges nor independent of the underlying graph. After a domination step (Lemma~\ref{lem:domination}), however, it can be replaced by an independent locally dependent percolation model in the following sense.

\paragraph{Locally dependent percolation models.}
Let $\mu$ be a probability measure on the subsets of $\mathcal N_o$. In the \emph{locally dependent percolation model} with local law $\mu$, every vertex $x$ samples a set $A_x\subset x+\mathcal N_o$, independently of all other vertices and such that $A_x-x$ has law $\mu$, and every undirected edge $\{x,y\}$ with $y\in A_x$ is declared open. We write $\mathbf P$ for the law of the resulting set $\omega\subset E(\Z^d)$ of open edges, where $E(\Z^d)$ denotes the set of undirected nearest-neighbor edges, and we say that $\mathbf P$ has \emph{$\delta$-marginals} if $\mu(\{A\colon y\in A\})\ge\delta$ for all $y\in\mathcal N_o$. For instance, the undirected $2dp$-neighbor graph is the locally dependent percolation model whose local law is the law of $A_o$ under $\P_p$; it has $p$-marginals. A set $X\subset E(\Z^d)$ is \emph{everywhere percolating} if every vertex lies in an infinite connected component of the graph $(\Z^d,X)$.

\begin{proposition}\label{proposition:locally-dependent-sprinkling}
Let $d\ge2$ and $\delta\in(0,1]$. There are constants $c_1,c_2>0$, depending only on $d$ and $\delta$, such that the following holds for every everywhere percolating set $X\subset E(\Z^d)$ and every locally dependent percolation model $\mathbf P$ with $\delta$-marginals. For all $n\in\N$, with $\mathbf P$-probability at least $1-c_1\exp(-c_2n)$, any two vertices of $\Lambda_n$ are joined by a path of edges of $X\cup\omega$ that stays inside $\Lambda_{2n}$.
\end{proposition}

The proof is inspired by the Markov chain argument of Contreras and Hirschi~\cite[Section~5.2]{contrerasThesis} for independent Bernoulli sprinkling. We reveal the sprinkling layer by layer, from $\partial\Lambda_{2n}$ inwards, and follow the number of components that meet the current box. Since $X$ percolates everywhere, every such component reaches the boundary of the current box. Along this boundary the components can be arranged in a tree, and each edge of the tree is opened with probability at least $\delta$ by a vertex whose choice has not been revealed yet, independently of the other edges of the tree. The number of components is therefore dominated by a coalescent in which each of the initial separations disappears independently with probability $\delta$ in every layer. Revealing the sprinkling vertex by vertex, rather than edge by edge, is what makes this work for dependent sprinkling.

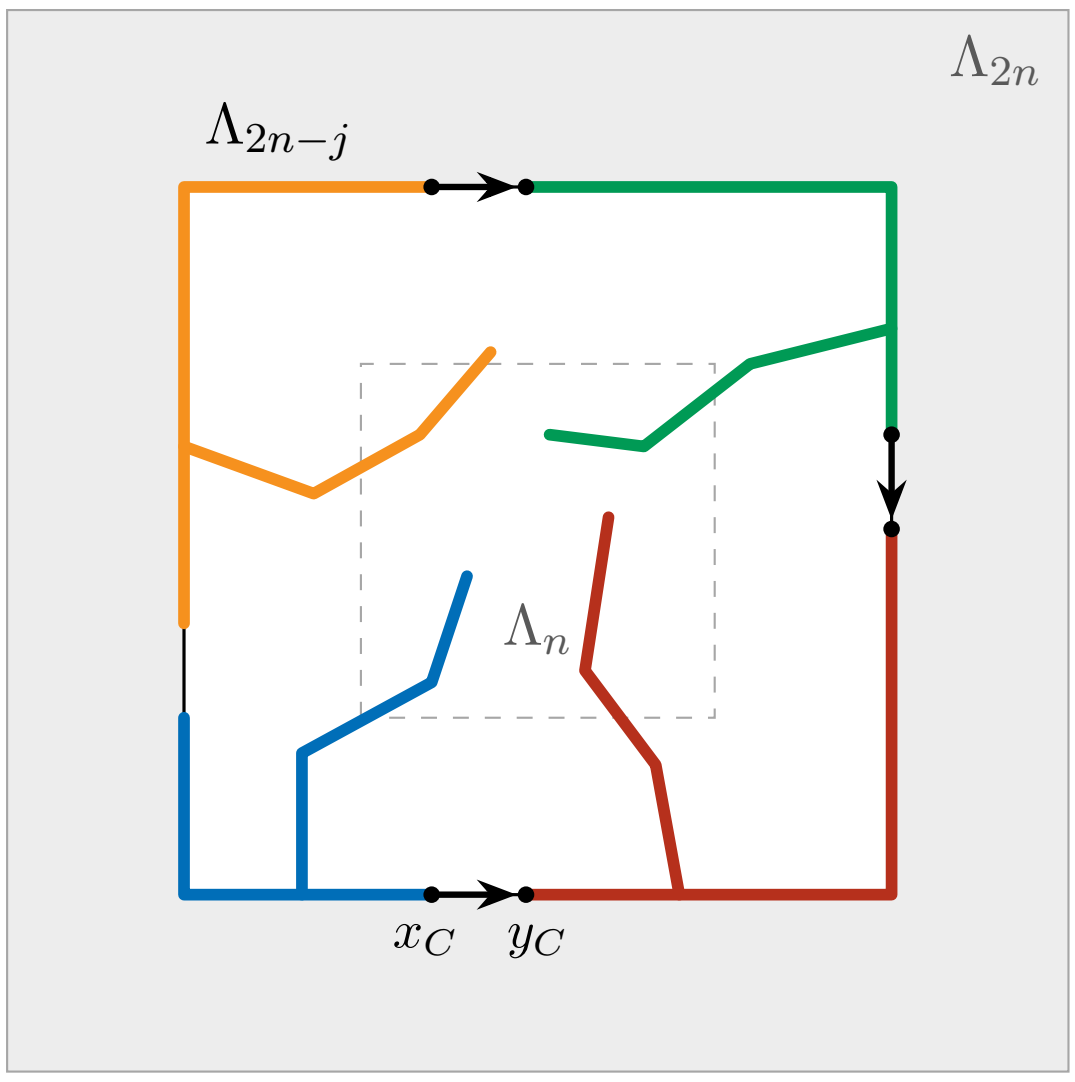
\begin{figure}[ht]
\centering
\begin{tikzpicture}[scale=0.78,
  comp/.style={line width=2.2pt, line cap=round, line join=round},
  bv/.style={circle, fill=black, inner sep=1.1pt},
  tr/.style={-{Stealth[length=2.6mm]}, line width=1.2pt}]
  \fill[gray!14, even odd rule] (-4.5,-4.5) rectangle (4.5,4.5) (-3,-3) rectangle (3,3);
  \draw[gray!70] (-4.5,-4.5) rectangle (4.5,4.5);
  \draw[black, line width=0.6pt] (-3,-3) rectangle (3,3);
  \draw[gray!70, dashed] (-1.5,-1.5) rectangle (1.5,1.5);
  \node[gray!70!black, anchor=north east] at (4.45,4.45) {$\Lambda_{2n}$};
  \node[anchor=south west] at (-3.0,3.05) {$\Lambda_{2n-j}$};
  \node[gray!70!black] at (0,-0.75) {$\Lambda_n$};
  \draw[comp, NavyBlue] (-3,-1.5) -- (-3,-3) -- (-0.9,-3);
  \draw[comp, NavyBlue] (-2,-3) -- (-2,-1.8) -- (-0.9,-1.2) -- (-0.6,-0.3);
  \draw[comp, BrickRed] (-0.1,-3) -- (3,-3) -- (3,0.1);
  \draw[comp, BrickRed] (1.2,-3) -- (1.0,-1.9) -- (0.4,-1.1) -- (0.6,0.2);
  \draw[comp, ForestGreen] (3,0.9) -- (3,3) -- (-0.1,3);
  \draw[comp, ForestGreen] (3,1.8) -- (1.8,1.5) -- (0.9,0.8) -- (0.1,0.9);
  \draw[comp, BurntOrange] (-0.9,3) -- (-3,3) -- (-3,-0.7);
  \draw[comp, BurntOrange] (-3,0.8) -- (-1.9,0.4) -- (-1.0,0.9) -- (-0.4,1.6);
  \draw[tr] (-0.9,-3) -- (-0.18,-3);   
  \draw[tr] (3,0.9) -- (3,0.18);       
  \draw[tr] (-0.9,3) -- (-0.18,3);     
  \foreach \p in {(-0.9,-3),(-0.1,-3),(3,0.9),(3,0.1),(-0.9,3),(-0.1,3)} \node[bv] at \p {};
  \node[below, font=\small] at (-0.95,-3.08) {$x_C$};
  \node[below, font=\small] at (0.0,-3.08) {$y_C$};
\end{tikzpicture}
\caption{One step of the sprinkling argument in the proof of Proposition~\ref{proposition:locally-dependent-sprinkling}. The edges selected by the vertices of the shaded annulus $\Lambda_{2n}\setminus\Lambda_{2n-j}$ have been revealed. Each component of $X\cup\omega_j$ that meets $\Lambda_{2n-j}$ reaches $\partial\Lambda_{2n-j}$; along $\partial\Lambda_{2n-j}$ the components are arranged in a tree, and each tree edge $\{x_C,y_C\}$ (arrows) is selected by the unrevealed vertex $x_C\in\partial\Lambda_{2n-j}$ with probability at least $\delta$, independently of the other tree edges.}\label{fig:sprinkling}
\end{figure}

\begin{proof}
Fix $n\in\N$, $X$ and $\mathbf P$. For $H\subset E(\Z^d)$ and $1\le m\le2n$, let $\mathcal C_m(H)$ be the set of connected components of the graph on $\Lambda_{2n}$ formed by the edges of $H$ with both endpoints in $\Lambda_{2n}$, that intersect $\Lambda_m$. Since $X$ is everywhere percolating, we have for every $H\supset X$
\begin{equation}\label{eq:reach-boundary}
\text{every component in }\mathcal C_m(H)\text{ contains a vertex of }\partial\Lambda_m\,.
\end{equation}
Indeed, a vertex $x\in\Lambda_m$ of the component lies in an infinite component of $(\Z^d,X)$, so a path of edges of $X$ leads from $x$ to a vertex outside $\Lambda_m$; up to its first vertex in $\partial\Lambda_m$ this path stays inside $\Lambda_m$, so it belongs to the component.

For $j\in\{0,\dots,n\}$ let $\omega_j$ be the set of edges $\{x,y\}$ with $x\in\Lambda_{2n}\setminus\Lambda_{2n-j}$ and $y\in A_x\cap\Lambda_{2n}$. These are the edges of $\omega$ inside $\Lambda_{2n}$ that are selected by the vertices of the annulus $\Lambda_{2n}\setminus\Lambda_{2n-j}$; in particular $\omega_0=\emptyset$. Let $\mathcal F_j$ be the $\sigma$-algebra generated by $(A_x)_{x\in\Lambda_{2n}\setminus\Lambda_{2n-j}}$, and let
\begin{equation}
K_j:=|\mathcal C_{2n-j}(X\cup\omega_j)|
\end{equation}
be the number of components of $X\cup\omega_j$ in $\Lambda_{2n}$ that intersect $\Lambda_{2n-j}$. Then $K_j$ is $\mathcal F_j$-measurable, $K_0\le|\partial\Lambda_{2n}|$ by~\eqref{eq:reach-boundary}, and if $K_n=1$, then any two vertices of $\Lambda_n$ are joined by a path of edges of $X\cup\omega_n\subset X\cup\omega$ inside $\Lambda_{2n}$.

\emph{Step 1: One layer.} We show that for $0\le j<n$, conditionally on $\mathcal F_j$, the random variable $K_{j+1}-1$ is stochastically dominated by a binomial random variable with parameters $K_j-1$ and $1-\delta$. Let $W:=\mathcal C_{2n-j}(X\cup\omega_j)$, so that $|W|=K_j$, and call two components in $W$ adjacent if they contain neighboring vertices of $\partial\Lambda_{2n-j}$. This graph on $W$ is connected: every component in $W$ contains a vertex of $\partial\Lambda_{2n-j}$ by~\eqref{eq:reach-boundary}, every vertex of $\partial\Lambda_{2n-j}$ belongs to a component in $W$, and $\partial\Lambda_{2n-j}$ is connected in $\Z^d$ since $d\ge2$. Choose, as a function of $X$ and $\omega_j$, a spanning tree of this graph and a root, and for every component $C\in W$ other than the root an edge $\{x_C,y_C\}$ with $x_C\in C\cap\partial\Lambda_{2n-j}$ and $y_C\in\partial\Lambda_{2n-j}$ in the parent of $C$; see Figure~\ref{fig:sprinkling}. The vertices $x_C$ are distinct, since they lie in distinct components, and they lie in $\partial\Lambda_{2n-j}=\Lambda_{2n-j}\setminus\Lambda_{2n-j-1}$. Hence their sets $A_{x_C}$ are independent of $\mathcal F_j$ and of each other, and $\{x_C,y_C\}\in\omega_{j+1}$ whenever $y_C\in A_{x_C}$. Let $Z$ be the number of non-root components $C$ with $y_C\in A_{x_C}$. Conditionally on $\mathcal F_j$, it is a sum of $K_j-1$ independent indicators with success probabilities at least $\delta$. The corresponding edges form a forest with $Z$ edges on $W$, so they merge the components in $W$ into at most $K_j-Z$ components. Every component of $X\cup\omega_{j+1}$ that intersects $\Lambda_{2n-j}$ contains a component in $W$, and therefore
\begin{equation}
K_{j+1}\le|\mathcal C_{2n-j}(X\cup\omega_{j+1})|\le K_j-Z\,.
\end{equation}
Since $K_j-1-Z$ is a sum of $K_j-1$ conditionally independent indicators with success probabilities at most $1-\delta$, Step 1 follows.

\emph{Step 2: Iteration.} For $\ell\in\N_0$ and $s\in[0,1]$ let $\mathrm{Bin}(\ell,s)$ denote a binomial random variable with parameters $\ell$ and $s$, defined on an extension of the probability space and independent of everything else, and let $\mathbf E$ denote the expectation on this space. We show by induction on $j$ that $K_j-1$ is stochastically dominated by $\mathrm{Bin}(K_0-1,(1-\delta)^j)$; note that $K_0$ is deterministic. For $j=0$ there is nothing to prove. For the induction step let $f$ be a nondecreasing function on $\N_0$ and set $g(\ell):=\mathbf E[f(\mathrm{Bin}(\ell,1-\delta))]$; then $g$ is nondecreasing. By Step~1 and the induction hypothesis,
\begin{equation}
\mathbf E\big[f(K_{j+1}-1)\big]\le\mathbf E\big[g(K_j-1)\big]\le\mathbf E\big[g\big(\mathrm{Bin}(K_0-1,(1-\delta)^j)\big)\big]=\mathbf E\big[f\big(\mathrm{Bin}(K_0-1,(1-\delta)^{j+1})\big)\big]\,,
\end{equation}
where in the last step we used that thinning a binomial random variable with parameters $K_0-1$ and $(1-\delta)^j$ independently with retention probability $1-\delta$ produces a binomial random variable with parameters $K_0-1$ and $(1-\delta)^{j+1}$.

\emph{Step 3: Conclusion.} By Step 2 with $j=n$, Bernoulli's inequality and $K_0\le|\partial\Lambda_{2n}|\le 2d(4n+1)^{d-1}$,
\begin{equation}
\mathbf P(K_n=1)\ge\big(1-(1-\delta)^n\big)^{K_0-1}\ge1-K_0(1-\delta)^n\ge1-2d(4n+1)^{d-1}(1-\delta)^n\,.
\end{equation}
The right-hand side is at least $1-c_1\exp(-c_2n)$ for suitable $c_1,c_2>0$ depending only on $d$ and $\delta$.
\end{proof}

To apply Proposition~\ref{proposition:locally-dependent-sprinkling}, we split the $2dp$-neighbor graph above the threshold into an everywhere percolating part and an independent sprinkling.
\begin{lemma}\label{lem:domination}
Let $d\ge2$ and $p>(d+1)/(2d)$, and let $q:=\min\{2dp-d-1,1\}/(2d)$. Let $X$ and $\omega$ be independent, where $X$ is the set of open edges of the undirected $(d+1)$-neighbor graph, that is, of the undirected $2dp$-neighbor graph with $p=(d+1)/(2d)$, and $\omega$ is the set of open edges of the undirected $2dq$-neighbor graph. Then $X\cup\omega$ is stochastically dominated by the set of open edges of the undirected $2dp$-neighbor graph.
\end{lemma}
\begin{proof}
Let $B$ and $B'$ be independent random subsets of $\mathcal N_o$ with the laws of $A_o$ under $\P_{(d+1,0)}$ and $\P_q$, respectively, and let $A$ have the law of $A_o$ under $\P_p$. The set $X\cup\omega$ is a function of the independent family of the sets selected by the vertices in either graph, which are distributed as translates of $B\cup B'$, since an edge of $X\cup\omega$ is selected by one of its endpoints in one of the two graphs. It therefore suffices to couple $B\cup B'$ and $A$ such that $B\cup B'\subset A$, and to use independent copies of this coupling at all vertices. The laws of $B$, $B'$ and $A$ are invariant under permutations of $\mathcal N_o$, and so is the law of $B\cup B'$. If two random subsets of $\mathcal N_o$ have permutation-invariant laws, and the cardinality of the first is stochastically dominated by the cardinality of the second, then the first set is stochastically dominated by the second: couple the cardinalities monotonically, choose the second set uniformly among the sets of its cardinality, and the first set uniformly among the subsets of the second set of its cardinality. It therefore suffices to compare cardinalities. Now $|B|=d+1$, and $B'$ is empty with probability $1-2dq$ and otherwise consists of a single uniformly chosen neighbor, so $|B\cup B'|\in\{d+1,d+2\}$ and
\begin{equation}
\P\big(|B\cup B'|=d+2\big)=2dq\cdot\frac{d-1}{2d}=(d-1)q\,.
\end{equation}
On the other hand, $|A|\ge d+1$ always, $|A|\ge d+2$ almost surely if $2dp\ge d+2$, and $\P(|A|=d+2)=2dp-d-1=2dq\ge(d-1)q$ if $d+1<2dp<d+2$. In both cases the cardinality of $B\cup B'$ is stochastically dominated by that of $A$.
\end{proof}

\begin{proof}[Proof of Theorem~\ref{thm:full-coverage-local-uniqueness}]
(i) Let $p\ge(d+1)/(2d)$. By Lemma~\ref{lem:forward-infinite}, almost surely every vertex lies in an infinite open cluster, and by Proposition~\ref{prop:BK} there is almost surely at most one infinite open cluster. Hence all vertices lie in the same open cluster.

(ii) Let $p>(d+1)/(2d)$ and let $X$, $\omega$ and $q>0$ be as in Lemma~\ref{lem:domination}. By Lemma~\ref{lem:forward-infinite}, applied with parameter $(d+1)/(2d)$, almost surely every vertex has an infinite forward cluster in the directed graph underlying $X$, so $X$ is almost surely everywhere percolating. Since $\omega$ is independent of $X$ and has $q$-marginals, Proposition~\ref{proposition:locally-dependent-sprinkling} with $\delta=q$, applied conditionally on $X$, shows that with probability at least $1-c_1\exp(-c_2n)$ any two vertices of $\Lambda_n$ are joined by a path of edges of $X\cup\omega$ inside $\Lambda_{2n}$. Here $c_1$ and $c_2$ depend only on $d$ and $q$, hence only on $d$ and $p$. This event is increasing in the set of open edges, so by Lemma~\ref{lem:domination} the same bound holds for the undirected $2dp$-neighbor graph.
\end{proof}

\begin{proof}[Proof of Corollary~\ref{cor:non-analytic}]
Let $p<(d+1)/(2d)$, so that $k\le d$. Let $I$ be the event that the out-set of every vertex of the unit cube $\{0,1\}^d$ is contained in the cube and that no vertex outside the cube selects a vertex of the cube. Every vertex of the cube has exactly $d$ neighbors in the cube. Its out-set has at most $d$ elements with positive probability (always if $k<d$, and with probability $1-\eps$ if $k=d$), and since out-sets of a given size are uniformly distributed, it is then contained in the cube with positive probability. Every vertex outside the cube has at most one neighbor in the cube, which it avoids with probability $1-p>0$. Since only finitely many vertices are involved, $\P_p(I)>0$. On $I$ the open cluster of the origin is contained in the cube, so $\theta^\textsc{U}_d(p)\le1-\P_p(I)<1$. Moreover, by ergodicity, almost surely the event $I$ occurs for some translate of the cube, and then the open clusters of the vertices of this translate are finite; hence $\calC^\textsc{U}(o)\ne\Z^d$ almost surely. For $p\ge(d+1)/(2d)$ the claims follow from Theorem~\ref{thm:full-coverage-local-uniqueness}(i) and Lemma~\ref{lem:forward-infinite}.

The function $\theta^\textsc{U}_d$ is continuous on $(p^\textsc{U}_{\rm c}(d),1]$ by Theorem~\ref{thm_con}(i). If it were real-analytic on $(p^\textsc{U}_{\rm c}(d),1)$, then, being equal to one on the interval $((d+1)/(2d),1)$, it would be equal to one on $(p^\textsc{U}_{\rm c}(d),1)$. This is impossible, since $\theta^\textsc{U}_d<1$ on $(p^\textsc{U}_{\rm c}(d),(d+1)/(2d))$, which is nonempty because $p^\textsc{U}_{\rm c}(d)\le p^\textsc{D}_{\rm c}(d)<1/2$ by~\cite{CHJK24}.
\end{proof}

\subsection{Directed and strong clusters in the plane}\label{sec:planar}
Throughout this section $d=2$. We use the planar crossing estimate of~\cite{CHJK24} to surround every box by a finite strongly connected set. Every infinite forward cluster starting inside the box reaches this set and hence contains all its vertices. This observation gives uniqueness and the stronger conclusions needed for continuity.

A finite set $S\subset\Z^2\setminus\Lambda_m$ \emph{separates $\Lambda_m$ from infinity} if only finitely many vertices can be reached from $\Lambda_m$ by a nearest-neighbor path avoiding $S$. Here these paths need not be open or follow edge orientations. A set $S$ is \emph{strongly connected} if $u\rightsquigarrow^\textsc{D}v$ for every $u,v\in S$.

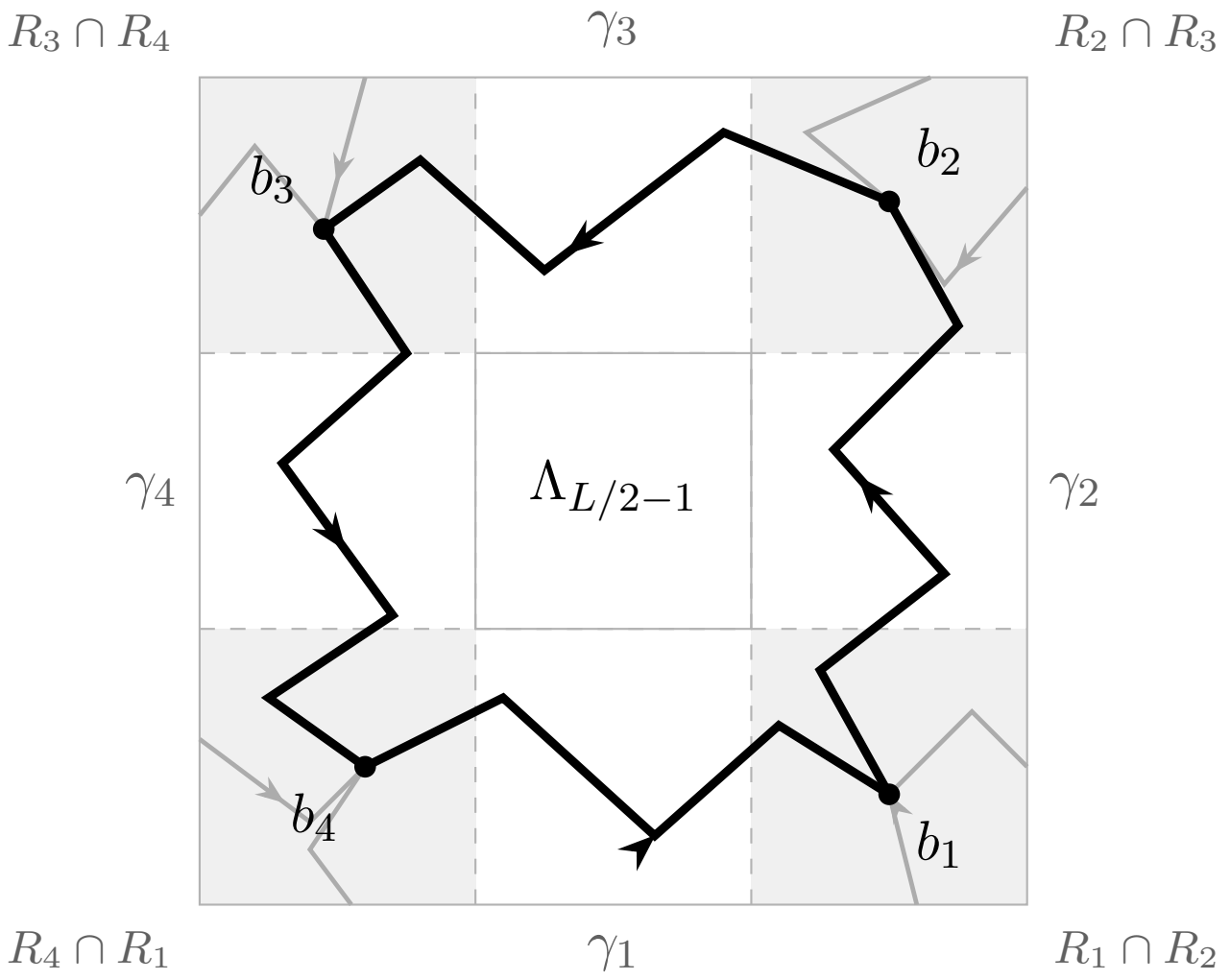
\begin{figure}[ht]
\centering
\begin{tikzpicture}[scale=0.95,
  cross/.style={gray!65, line width=0.9pt},
  crossarrow/.style={postaction={decorate},decoration={markings,mark=at position 0.1 with {\arrow[gray!65]{Stealth[length=2mm]}}}},
  walk/.style={black, line width=1.7pt},
  mid arrow/.style={postaction={decorate},decoration={markings,mark=at position 0.55 with {\arrow{Stealth[length=2.6mm]}}}},
  pt/.style={circle, fill=black, inner sep=1.5pt}]
  \fill[gray!12] (-3,-3) rectangle (-1,-1);
  \fill[gray!12] (1,-3) rectangle (3,-1);
  \fill[gray!12] (1,1) rectangle (3,3);
  \fill[gray!12] (-3,1) rectangle (-1,3);
  \draw[gray!60] (-3,-3) rectangle (3,3);
  \draw[gray!60] (-1,-1) rectangle (1,1);
  \draw[gray!60, dashed] (-3,-1) -- (3,-1);
  \draw[gray!60, dashed] (-3,1) -- (3,1);
  \draw[gray!60, dashed] (-1,-3) -- (-1,3);
  \draw[gray!60, dashed] (1,-3) -- (1,3);
  \node at (0,0) {$\Lambda_{L/2-1}$};
  \coordinate (b1) at (2.0,-2.2);
  \coordinate (b2) at (2.0,2.1);
  \coordinate (b3) at (-2.1,1.9);
  \coordinate (b4) at (-1.8,-2.0);
  \draw[cross, crossarrow] (-3,-1.8) -- (-2.2,-2.4) -- (b4) -- (-0.8,-1.5) -- (0.3,-2.5) -- (1.2,-1.7) -- (b1) -- (2.6,-1.6) -- (3,-2.0);
  \draw[cross, crossarrow] (2.2,-3) -- (b1) -- (1.5,-1.3) -- (2.4,-0.6) -- (1.6,0.3) -- (2.5,1.2) -- (b2) -- (1.4,2.6) -- (2.3,3);
  \draw[cross, crossarrow] (3,2.2) -- (2.4,1.5) -- (b2) -- (0.8,2.6) -- (-0.5,1.6) -- (-1.4,2.4) -- (b3) -- (-2.6,2.5) -- (-3,2.0);
  \draw[cross, crossarrow] (-1.8,3) -- (b3) -- (-1.5,1.0) -- (-2.4,0.2) -- (-1.6,-0.9) -- (-2.5,-1.5) -- (b4) -- (-2.2,-2.6) -- (-1.9,-3);
  \draw[walk, mid arrow] (b4) -- (-0.8,-1.5) -- (0.3,-2.5) -- (1.2,-1.7) -- (b1);
  \draw[walk, mid arrow] (b1) -- (1.5,-1.3) -- (2.4,-0.6) -- (1.6,0.3) -- (2.5,1.2) -- (b2);
  \draw[walk, mid arrow] (b2) -- (0.8,2.6) -- (-0.5,1.6) -- (-1.4,2.4) -- (b3);
  \draw[walk, mid arrow] (b3) -- (-1.5,1.0) -- (-2.4,0.2) -- (-1.6,-0.9) -- (-2.5,-1.5) -- (b4);
  \foreach \b in {b1,b2,b3,b4} \node[pt] at (\b) {};
  \node[anchor=north west] at ($(b1)+(0.05,-0.05)$) {$b_1$};
  \node[anchor=south west] at ($(b2)+(0.05,0.05)$) {$b_2$};
  \node[anchor=south east] at ($(b3)+(-0.05,0.05)$) {$b_3$};
  \node[anchor=north east] at ($(b4)+(-0.05,-0.05)$) {$b_4$};
  \node[gray!75!black] at (0.0,-3.35) {$\gamma_1$};
  \node[gray!75!black] at (3.35,0.0) {$\gamma_2$};
  \node[gray!75!black] at (0.0,3.35) {$\gamma_3$};
  \node[gray!75!black] at (-3.35,0.0) {$\gamma_4$};
  \node[gray!75!black, font=\footnotesize, anchor=north east] at (-3.05,-3.05) {$R_4\cap R_1$};
  \node[gray!75!black, font=\footnotesize, anchor=north west] at (3.05,-3.05) {$R_1\cap R_2$};
  \node[gray!75!black, font=\footnotesize, anchor=south west] at (3.05,3.05) {$R_2\cap R_3$};
  \node[gray!75!black, font=\footnotesize, anchor=south east] at (-3.05,3.05) {$R_3\cap R_4$};
\end{tikzpicture}
\caption{The four directed crossings $\gamma_1,\dots,\gamma_4$, oriented counterclockwise around $\Lambda_{L/2-1}$ (schematic). Joining them at the corner vertices $b_i$ gives a closed directed path (bold). Shortening this walk as in the proof of Lemma~\ref{lem:gluing} gives a simple directed cycle separating the box from infinity.}\label{fig:gluing}
\end{figure}

\paragraph{Surrounding a box.}
A \emph{left-right crossing} of a rectangle is an open directed path contained in the rectangle that joins its left side to its right side; the other crossing directions are defined analogously. The proof of~\cite[Lemma~3.2]{CHJK24} gives
\begin{equation}\label{eq:crossing}
\lim_{L\to\infty}\P_{1/2}\big(\text{there is a left-right crossing of }[0,3L]\times[0,L]\big)=1\,.
\end{equation}
Four suitably oriented crossings can be joined around a box, as in Figure~\ref{fig:gluing}. Let us make precise what we mean by that. 

\begin{lemma}\label{lem:gluing}
Let $L\ge4$ be even and $\omega\in\Omega$. Set $R_1:=[-3L/2,3L/2]\times[-3L/2,-L/2]$, and let $R_2,R_3,R_4$ be its rotations about the origin by $\pi/2$, $\pi$ and $3\pi/2$. Suppose these rectangles have open directed crossings from left to right, bottom to top, right to left and top to bottom, respectively. Then there is a closed directed walk of open edges whose vertex set is contained in $R_1\cup\dots\cup R_4$ and separates $\Lambda_{L/2-1}$ from infinity.
\end{lemma}

We believe a picture like in Figure ~\ref{fig:gluing} explains the idea better than words can and therefore omit the proof of this lemma. 

\begin{proposition}\label{prop:cont-planar}
Let $d=2$ and $p\in[1/2,1]$. Then $\P_p$-almost surely, for every $m\in\N$ there is a finite set $S\subset\Z^2\setminus\Lambda_m$ that is strongly connected and separates $\Lambda_m$ from infinity.
\end{proposition}
\begin{proof}
Fix $m$ and let $L\ge2m+2$ be even. By~\eqref{eq:crossing}, translation and rotation invariance, and the monotone coupling of Lemma~\ref{lem:cont-coupling}, each crossing required in Lemma~\ref{lem:gluing} has $\P_p$-probability tending to one as $L\to\infty$. A union bound and that lemma give, with probability tending to one, a closed directed walk outside $\Lambda_m$ whose vertex set separates $\Lambda_m$ from infinity. Its vertex set is strongly connected. Thus a suitable set $S$ exists almost surely for this $m$; intersecting over $m\in\N$ proves the claim.
\end{proof}

\paragraph{Consequences for infinite clusters.}
The following deterministic lemma gives uniqueness and the three inputs used in the continuity proof. In particular, part~(i) says that any two infinite forward clusters differ only in finitely many vertices.
\begin{lemma}\label{lem:cont-planar}
Let $d=2$ and let $\omega\in\Omega$ be such that for every $m\in\N$ there is a finite set $S_m\subset\Z^2\setminus\Lambda_m$ which is strongly connected in $\omega$ and separates $\Lambda_m$ from infinity. Then the following hold in $\omega$ for all $x,y\in\Z^2$.
\begin{enumerate}[(i)]
\item If $|\calC^\textsc{D}(x)|=|\calC^\textsc{D}(y)|=\infty$, then $\calC^\textsc{D}(x)\setminus\calC^\textsc{D}(y)$ is finite.
\item If $|\calC^\textsc{S}(x)|=|\calC^\textsc{S}(y)|=\infty$, then $\calC^\textsc{S}(x)=\calC^\textsc{S}(y)$.
\item $|\calC^\textsc{S}(x)|=\infty$ if and only if $|\calC^\textsc{D}(x)|=\infty$ and $|\calC^{\leftarrow}(x)|=\infty$.
\end{enumerate}
In particular, $N^\textsc{D}\le1$ and $N^\textsc{S}\le1$ in $\omega$.
\end{lemma}
\begin{proof}
Fix $m$ and let $D_m$ be the finite set of vertices reachable from $\Lambda_m$ by a nearest-neighbor path avoiding $S_m$. If $x\in\Lambda_m$ has an infinite forward cluster, an open directed path from $x$ to a vertex outside $D_m$ meets $S_m$. Strong connectivity then implies that $x$ reaches every vertex of $S_m$. Applying the same argument with all arrows reversed shows that every vertex of $S_m$ reaches $x$ if its backward cluster is infinite. In particular, if both clusters are infinite, then $S_m\subset\calC^\textsc{S}(x)$.

For (i), choose $m$ with $x,y\in\Lambda_m$. Every open directed path from $x$ to a vertex $u\notin D_m$ meets $S_m$. Since $y$ reaches every vertex of $S_m$, every such $u\in\calC^\textsc{D}(x)$ also belongs to $\calC^\textsc{D}(y)$. Hence $\calC^\textsc{D}(x)\setminus\calC^\textsc{D}(y)\subset D_m$ is finite.

For (ii), choose $m$ with $x,y\in\Lambda_m$. Infinite strong forward clusters have infinite forward and backward clusters, so both $\calC^\textsc{S}(x)$ and $\calC^\textsc{S}(y)$ contain the nonempty set $S_m$. As strong forward clusters are equivalence classes, they coincide.

For (iii), one implication follows from $\calC^\textsc{S}(x)=\calC^\textsc{D}(x)\cap\calC^{\leftarrow}(x)$. Conversely, if both clusters are infinite, the initial observation gives $S_m\subset\calC^\textsc{S}(x)$ for every sufficiently large $m$. Since $S_m$ is nonempty and lies outside $\Lambda_m$, the strong forward cluster is infinite. The assertion about $N^\textsc{D}$ follows from (i), applied also with the roles of $x$ and $y$ exchanged, and the assertion about $N^\textsc{S}$ follows from (ii).
\end{proof}

\begin{proof}[Proof of Theorem~\ref{thm_uni}(ii)]
For $p<p^*_{\rm c}(2)$, almost surely every cluster $\calC^*(x)$ is finite, so $N^*=0$. For $p\ge1/2$, Proposition~\ref{prop:cont-planar} and Lemma~\ref{lem:cont-planar} give $N^\textsc{D}\le1$ and $N^\textsc{S}\le1$ almost surely.
\end{proof}
\subsection{Continuity of the percolation functions}\label{sec:continuity}
We follow the classical route of van den Berg and Keane~\cite{vdBK84} and of Aizenman, Kesten and Newman~\cite{aizenman1987uniqueness}, see also~\cite[Chapter~8]{Grimmett1999}. For forward and undirected clusters, right-continuity holds because percolation is a decreasing limit of events that depend on finitely many edges; for strong forward clusters we combine this with the planar structure of Section~\ref{sec:planar}. Left-continuity at a supercritical parameter $p$ rests on uniqueness: a finite path that is open at level $p$ is open at a slightly smaller level, and uniqueness at level $p$ connects the origin by such a path to an infinite cluster of a lower level. Both steps take place in the coupling of Lemma~\ref{lem:cont-coupling}.

Two consequences of this coupling are used repeatedly. First, since $\mathcal E$ is countable and the labels have no atoms by Lemma~\ref{lem:cont-coupling}(iii), for every fixed $p\in[0,1]$ almost surely no label equals $p$. Therefore, for every $p\in[0,1)$,
\begin{equation}\label{eq:cont-pplus}
\lim_{p'\downarrow p}\omega_{p'}=\big(\one\{T_e\le p\}\big)_{e\in\mathcal E}=\omega_p\qquad\P\text{-almost surely}\,.
\end{equation}
Second, if a finite directed path with at least one edge is open at level $p$, that is, all its edges are open in $\omega_p$, then the largest label on the path is smaller than $p$, and the path is open at every level above this largest label, in particular at some level smaller than $p$; a path without edges is open at every level.

\begin{lemma}\label{lem:cont-right}
Let $(E_m)_{m\ge1}$ be increasing events such that $E_m$ depends only on the edges with both endpoints in $\Lambda_m$, and let $E:=\bigcap_{m\ge1}E_m$. Then for every $p\in[0,1)$ and every sequence $p_n\downarrow p$ in $(p,1]$,
\begin{equation}\label{eq:cont-right}
\bigcap_{n\ge1}\{\omega_{p_n}\in E\}=\{\omega_p\in E\}\qquad\P\text{-almost surely}\,,
\end{equation}
and in particular $p\mapsto\P(\omega_p\in E)$ is right-continuous on $[0,1)$. This applies to the three events $\{|\calC^\textsc{D}(o)|=\infty\}$, $\{|\calC^\textsc{U}(o)|=\infty\}$ and $\{|\calC^\textsc{D}(o)|=|\calC^{\leftarrow}(o)|=\infty\}$.
\end{lemma}
\begin{proof}
The inclusion $\supset$ in~\eqref{eq:cont-right} holds since $\omega_p\le\omega_{p_n}$ and $E$ is increasing. For the inclusion $\subset$, assume that $\omega_{p_n}\in E$ for all $n$ and fix $m$. For each of the finitely many edges $e$ with both endpoints in $\Lambda_m$, the sequence $\omega_{p_n}(e)=\one\{T_e<p_n\}$ is nonincreasing and $\{0,1\}$-valued, hence eventually equal to $\one\{T_e\le p\}$. Since $\omega_{p_n}\in E_m$ for all $n$ and $E_m$ depends only on these edges, the configuration $(\one\{T_e\le p\})_{e\in\mathcal E}$ lies in $E_m$. As $m$ was arbitrary, it lies in $E$, and by~\eqref{eq:cont-pplus} it coincides with $\omega_p$ outside a null event that does not depend on the sequence $(p_n)_{n\ge1}$. The events $\{\omega_{p_n}\in E\}$ decrease in $n$, so~\eqref{eq:cont-right} gives $\P(\omega_{p_n}\in E)\downarrow\P(\omega_p\in E)$, which is the right-continuity.

For the three events, let $E_m:=\{o\rightsquigarrow^\textsc{D}\partial\Lambda_m\}$ be the event that some open directed path leads from $o$ to a vertex of $\partial\Lambda_m$. Cutting such a path at its first vertex in $\partial\Lambda_m$ leaves an open directed path all of whose vertices lie in $\Lambda_m$, so $E_m$ depends only on the edges with both endpoints in $\Lambda_m$, and $E_m$ is increasing. If $\calC^\textsc{D}(o)$ is infinite, it contains a vertex outside $\Lambda_{m-1}$, and an open directed path from $o$ to this vertex meets $\partial\Lambda_m$ when it first leaves $\Lambda_{m-1}$; conversely, if $E_m$ occurs for every $m$, then $\calC^\textsc{D}(o)$ meets each of the disjoint sets $\partial\Lambda_m$. Hence $\{|\calC^\textsc{D}(o)|=\infty\}=\bigcap_mE_m$. The backward cluster is treated in the same way with the events $\{\partial\Lambda_m\rightsquigarrow^\textsc{D}o\}$, cutting paths at their last vertex in $\partial\Lambda_m$, and intersecting the two families of events handles the third event. The undirected cluster is treated with the events that $o$ is joined to $\partial\Lambda_m$ by an undirected open path, which again depend only on the directed edges with both endpoints in $\Lambda_m$.
\end{proof}

\begin{proof}[Proof of Theorem~\ref{thm_con}(i)]
Given Theorem~\ref{thm_uni}(i), the argument is the classical one for Bernoulli percolation~\cite{vdBK84,aizenman1987uniqueness}, and we only indicate the input from our model. In the coupling of Lemma~\ref{lem:cont-coupling}, the undirected edge $\{x,y\}$ is open in $\omega_p$ if and only if $\min\{T_{(x,y)},T_{(y,x)}\}<p$, and this minimum of two independent uniform random variables has no atoms; so the coupling has the two properties stated before Lemma~\ref{lem:cont-right} also for undirected paths. Right-continuity of $\theta^\textsc{U}_d$ on $[0,1)$ is Lemma~\ref{lem:cont-right}. Left-continuity at $p\in(p^\textsc{U}_{\rm c}(d),1]$ follows exactly as in~\cite{vdBK84}: if the origin percolates at level $p$, then by Theorem~\ref{thm_uni}(i) at level $p$ its open cluster contains the infinite cluster of any lower supercritical level $p_1<p$, which exists by ergodicity so the origin is joined to it by a finite path that is open at some level below $p$. Finally, $\theta^\textsc{U}_d=0$ on $[0,p^\textsc{U}_{\rm c}(d))$.
\end{proof}

In the directed case, the left-continuity argument uses the bound $N^\textsc{D}\le1$ of Theorem~\ref{thm_uni}(ii) in the form of Lemma~\ref{lem:cont-planar}(i): two infinite forward clusters differ only in finitely many vertices, so the forward cluster of the origin at level $p$ reaches a vertex from which a lower level percolates. The weaker statement that any two infinite forward clusters intersect would not suffice (Remark~\ref{rem:cont-ND}).

\begin{proof}[Proof of Theorem~\ref{thm_con}(ii)]
Let $d=2$. By Proposition~\ref{prop:cont-planar}, for every $p\in[1/2,1]$ the configuration $\omega_p$ almost surely satisfies the assumption of Lemma~\ref{lem:cont-planar}; whenever we apply parts (i)--(iii) of that lemma at such a level $p$, we tacitly discard the corresponding null event. The lemma is only applied at the level $p$ at which continuity is proved, never at the approximating levels $p_n$, which may be smaller than $1/2$. Since $\theta^*_2=0$ on $[0,p^*_{\rm c}(2))$ and $\theta^*_2$ is monotone, it suffices to prove right-continuity at every $p\in[1/2,1)$ and left-continuity at every $p\in[1/2,1]$ with $p>p^*_{\rm c}(2)$, along one sequence each.

\emph{Case $*=\textsc{D}$, right-continuity.} For $p\in[0,1)$, Lemma~\ref{lem:cont-right} applied to $\{|\calC^\textsc{D}(o)|=\infty\}$ gives $\lim_{p'\downarrow p}\theta^\textsc{D}_2(p')=\theta^\textsc{D}_2(p)$.

\emph{Case $*=\textsc{D}$, left-continuity.} Let $p\in[1/2,1]$; then $p>p^\textsc{D}_{\rm c}(2)$. Choose $p_1<p_2<\dots$ in $(p^\textsc{D}_{\rm c}(2),p)$ with $p_n\uparrow p$. The events $\{|\calC^\textsc{D}_{p_n}(o)|=\infty\}$ increase in $n$, so $\lim_n\theta^\textsc{D}_2(p_n)=\P\big(\bigcup_n\{|\calC^\textsc{D}_{p_n}(o)|=\infty\}\big)$, and since $\theta^\textsc{D}_2(p_n)\le\theta^\textsc{D}_2(p)$ it suffices to show that almost surely
\begin{equation}\label{eq:cont-D-incl}
\{|\calC^\textsc{D}_p(o)|=\infty\}\subset\bigcup_{n\ge1}\{|\calC^\textsc{D}_{p_n}(o)|=\infty\}\,.
\end{equation}
By ergodicity there is almost surely a vertex $z$ with $|\calC^\textsc{D}_{p_1}(z)|=\infty$, and by Lemma~\ref{lem:cont-koenig} there is an infinite self-avoiding directed path $\rho$ starting at $z$ which is open at level $p_1$. Assume that $|\calC^\textsc{D}_p(o)|=\infty$. The vertices of $\rho$ lie in $\calC^\textsc{D}_{p_1}(z)\subset\calC^\textsc{D}_p(z)$, and $\calC^\textsc{D}_p(z)\setminus\calC^\textsc{D}_p(o)$ is finite by Lemma~\ref{lem:cont-planar}(i) at level $p$. Hence some vertex $w$ of $\rho$ lies in $\calC^\textsc{D}_p(o)$. Let $\gamma$ be a directed path from $o$ to $w$ which is open at level $p$, and choose $n$ such that $p_n$ exceeds all labels on $\gamma$. Then $\gamma$ is open at level $p_n$, and so is the infinite part of $\rho$ after $w$, which is open already at level $p_1\le p_n$. Therefore $|\calC^\textsc{D}_{p_n}(o)|=\infty$, which proves~\eqref{eq:cont-D-incl}.

\emph{Case $*=\textsc{S}$, left-continuity.} Let $p\in[1/2,1]$ with $p>p^\textsc{S}_{\rm c}(2)$, and choose $p_1<p_2<\dots$ in $(p^\textsc{S}_{\rm c}(2),p)$ with $p_n\uparrow p$. As in the directed case, it suffices to show that almost surely $\{|\calC^\textsc{S}_p(o)|=\infty\}\subset\bigcup_n\{|\calC^\textsc{S}_{p_n}(o)|=\infty\}$. By ergodicity there is almost surely a vertex $z$ with $|\calC^\textsc{S}_{p_1}(z)|=\infty$. Assume that $|\calC^\textsc{S}_p(o)|=\infty$. Since $\calC^\textsc{S}_{p_1}(z)\subset\calC^\textsc{S}_p(z)$, Lemma~\ref{lem:cont-planar}(ii) at level $p$ gives $\calC^\textsc{S}_p(o)=\calC^\textsc{S}_p(z)$, so there are directed paths $\gamma_1$ from $o$ to $z$ and $\gamma_2$ from $z$ to $o$ which are open at level $p$. Choose $n$ such that $p_n$ exceeds all labels on $\gamma_1$ and $\gamma_2$. Then $o$ and $z$ are joined in both directions by paths that are open at level $p_n$, so $o\in\calC^\textsc{S}_{p_n}(z)\supset\calC^\textsc{S}_{p_1}(z)$, which is infinite.

\emph{Case $*=\textsc{S}$, right-continuity.} Let $p\in[1/2,1)$ and $p_n\downarrow p$ in $(p,1]$. The events $\{|\calC^\textsc{S}_{p_n}(o)|=\infty\}$ decrease in $n$, so $\lim_n\theta^\textsc{S}_2(p_n)=\P(E')$ with $E':=\bigcap_n\{|\calC^\textsc{S}_{p_n}(o)|=\infty\}$, and $\theta^\textsc{S}_2(p)\le\P(E')$ by monotonicity. For the reverse inequality note that $\{|\calC^\textsc{S}_{p_n}(o)|=\infty\}\subset\{|\calC^\textsc{D}_{p_n}(o)|=|\calC^{\leftarrow}_{p_n}(o)|=\infty\}$ for every $n$, so that by Lemma~\ref{lem:cont-right}, applied to the third of the events listed there,
\begin{equation}
E'\subset\bigcap_{n\ge1}\big\{|\calC^\textsc{D}_{p_n}(o)|=|\calC^{\leftarrow}_{p_n}(o)|=\infty\big\}=\big\{|\calC^\textsc{D}_p(o)|=|\calC^{\leftarrow}_p(o)|=\infty\big\}\qquad\P\text{-almost surely}\,.
\end{equation}
By Lemma~\ref{lem:cont-planar}(iii) at level $p$, the event on the right-hand side almost surely coincides with $\{|\calC^\textsc{S}_p(o)|=\infty\}$. Hence $\P(E')\le\theta^\textsc{S}_2(p)$.
\end{proof}

\begin{remark}\label{rem:cont-ND}
The statement that any two infinite forward clusters intersect is weaker than $N^\textsc{D}\le1$ and does not suffice for the left-continuity argument. Consider, pathwise, the configuration in $\Z^2$ whose open edges are those of the two infinite directed paths $\pi=(\pi_0,\pi_1,\dots)$ and $\rho=(\rho_0,\rho_1,\dots)$ with $\pi_j:=(j,1)$ and $\rho_j:=(j,0)$, together with the edges $(\rho_j,\pi_j)$, $j\ge0$. Any two infinite forward clusters of this configuration intersect, but no vertex of $\pi$ has a directed path to $\rho$. If the labels of the edges of $\pi$ increase to $p$ along $\pi$ while $\rho$ is open at some level $p_1<p$, then $\pi_0$ percolates at level $p$ but at no smaller level. Here $N^\textsc{D}=2$, since $\calC^\textsc{D}(\rho_0)\setminus\calC^\textsc{D}(\pi_0)$ is infinite, and Lemma~\ref{lem:cont-planar}(i) excludes this configuration.
\end{remark}

\begin{remark}\label{rem:cont-critical}
For $*\in\{\textsc{D},\textsc{U}\}$ and every $d\ge2$, Lemma~\ref{lem:cont-right} gives right-continuity of $\theta^*_d$ at every $p\in[0,1)$, including $p=p^*_{\rm c}(d)$. Hence $\theta^*_d$ is continuous at $p^*_{\rm c}(d)$ if and only if $\theta^*_d(p^*_{\rm c}(d))=0$, exactly as for Bernoulli percolation. In the plane, the same holds for $\theta^\textsc{S}_2$ provided that $p^\textsc{S}_{\rm c}(2)\ge1/2$.
\end{remark}

\subsection{The order of the full-coverage transition}\label{sec:order}
In this section we prove Theorem~\ref{thm:order}. Throughout, $d\ge2$, and until the proof of Theorem~\ref{thm:order} at the end of the section, $p\in[1/2,(d+1)/(2d)]$. For these parameters every out-set has $d$ or $d+1$ elements, namely
\begin{equation}\label{eq:order-sizes}
\P_p(|A_x|=d)=d+1-2dp\qquad\text{and}\qquad\P_p(|A_x|=d+1)=2dp-d\,,
\end{equation}
since $\P_p=\P_{(d,2dp-d)}$ for $p<(d+1)/(2d)$ and $\P_p=\P_{(d+1,0)}$ for $p=(d+1)/(2d)$. As out-sets of a given size are uniformly distributed, a prescribed set of $d$ neighbors of $x$ is the out-set of $x$ with probability $(d+1-2dp)/\binom{2d}{d}$, and a prescribed set of $d+1$ neighbors with probability $(2dp-d)/\binom{2d}{d+1}$. Both are affine functions of $p$.

By Lemma~\ref{lem:forward-infinite}, $\theta^\textsc{U}_d=1$ on $[(d+1)/(2d),1]$, so Theorem~\ref{thm:order} concerns the probability $1-\theta^\textsc{U}_d(p)=\P_p(|\calC^\textsc{U}(o)|<\infty)$ for $p\in[1/2,(d+1)/(2d)]$, and in particular how it vanishes as $p\uparrow(d+1)/(2d)$. 

The key observation is the following. Every vertex $x$ of a finite open cluster $C$ satisfies $A_x\subset C$, so it has at least $d$ neighbors in $C$. Lemma~\ref{lem:corners} below shows that each finite set $C$ with this property has at least $2^d$ corners, i.e., vertices that lie in $C$ with exactly $d$ neighbors in $C$. In particular, each of these corners has to select exactly these $d$ neighbors to disconnect itself from the outside.

To exploit this, we decompose the event $\{|\calC^\textsc{U}(o)|<\infty\}$ according to the smallest box $\Lambda_n$ that contains $\calC^\textsc{U}(o)$. The probability of each part is a polynomial in $p$ whose degree grows only polynomially in $n$, and it is divisible by $(d+1-2dp)^{2^d}$, since the probability of any configuration of the out-sets in $\Lambda_{n+1}$ that confines the open cluster of the origin to $\Lambda_n$ contains this factor. A Peierls argument (Proposition~\ref{prop:tail}) shows that, after division by $(d+1-2dp)^{2^d}$, this polynomial is exponentially small in $n$, uniformly in $p$. Markov's inequality for derivatives of polynomials extends this bound to derivatives of every order, at a cost that is polynomial in $n$, so the series can be differentiated term by term, yielding the factorisation claimed in Proposition~\ref{prop:factorization}.

\begin{lemma}\label{lem:corners}
Let $S\subset\Z^d$ be finite and nonempty, and assume that every vertex of $S$ has at least $d$ neighbors in $S$. Then at least $2^d$ vertices of $S$ have exactly $d$ neighbors in $S$.
\end{lemma}
\begin{proof}
Let $e_1,\dots,e_d$ denote the unit vectors. For $s=(s_1,\dots,s_d)\in\{-1,1\}^d$ let $v_s$ be a vertex of $S$ at which the function $x\mapsto s_1x_1+\dots+s_dx_d$ attains its maximum over $x=(x_1,\dots,x_d)\in S$. As in the proof of Lemma~\ref{lem:forward-infinite}, none of the vertices $v_s+s_ie_i$, $i\in\{1,\dots,d\}$, lies in $S$, so the neighbors of $v_s$ in $S$ are among the $d$ vertices $v_s-s_ie_i$. By assumption, all of them lie in $S$, and $v_s$ has exactly $d$ neighbors in $S$. If $s\neq t$, say $t_i=-s_i$, then $v_s-s_ie_i\in S$, whereas $v_t-s_ie_i=v_t+t_ie_i\notin S$; hence $v_s\neq v_t$. Thus the $2^d$ vertices $v_s$ are distinct.
\end{proof}

\paragraph{Cost of closed boundaries.}
Recall that $E(\Z^d)$ denotes the set of undirected nearest-neighbor edges. For $S\subset\Z^d$ let $\partial_ES:=\{\{x,y\}\in E(\Z^d)\colon x\in S,\ y\notin S\}$ denote its edge boundary. We say that $\partial_ES$ is \emph{closed} if none of its edges is open in the undirected $2dp$-neighbor graph, that is, if $y\notin A_x$ and $x\notin A_y$ whenever $x\in S$, $y\notin S$ and $\{x,y\}\in\partial_ES$. If $\calC^\textsc{U}(o)$ is finite, then $\partial_E\calC^\textsc{U}(o)$ is closed.

\begin{lemma}\label{lem:cut}
For every finite nonempty set $S\subset\Z^d$,
\begin{equation}\label{eq:cut}
\P_p(\partial_ES\text{ is closed})\le\bigg(2^d(d+1-2dp)\binom{2d}{d}^{-1}\bigg)^{2^d}4^{-|\partial_ES|}\,.
\end{equation}
\end{lemma}
\begin{proof}
For $x\in S$ let $Y_x:=(x+\mathcal N_o)\setminus S$, and for $x\notin S$ let $Y_x:=(x+\mathcal N_o)\cap S$; thus $Y_x$ is the set of neighbors of $x$ on the other side of the boundary. Then $\partial_ES$ is closed if and only if $A_x\cap Y_x=\emptyset$ for all $x\in\Z^d$, and since every edge of $\partial_ES$ is counted once from each of its endpoints,
\begin{equation}\label{eq:cut-count}
\sum_{x\in\Z^d}|Y_x|=2|\partial_ES|\,.
\end{equation}
Let $x\in\Z^d$ and $t:=|Y_x|$. For $j\in\{d,d+1\}$, a uniformly chosen set of $j$ neighbors of $x$ avoids $Y_x$ with probability
\begin{equation}
\binom{2d-t}{j}\bigg/\binom{2d}{j}=\prod_{i=0}^{t-1}\frac{2d-j-i}{2d-i}\le\prod_{i=0}^{t-1}\frac{d-i}{2d-i}\le2^{-t}\,,
\end{equation}
and therefore $\P_p(A_x\cap Y_x=\emptyset)\le2^{-t}$. If $t>d$, this probability vanishes. If $t=d$, the only possible out-set of $x$ that avoids $Y_x$ is the set of the $d$ neighbors of $x$ outside $Y_x$, so that $\P_p(A_x\cap Y_x=\emptyset)=(d+1-2dp)/\binom{2d}{d}$.

If some vertex of $S$ has fewer than $d$ neighbors in $S$, then $|Y_x|>d$ for this vertex $x$, and~\eqref{eq:cut} holds trivially. Otherwise, by Lemma~\ref{lem:corners} there is a set $V\subset S$ of $2^d$ vertices with exactly $d$ neighbors in $S$, that is, with $|Y_x|=d$ for $x\in V$. Since the out-sets of different vertices are independent, the bounds above and~\eqref{eq:cut-count} give
\begin{equation}
\P_p(\partial_ES\text{ is closed})\le\bigg(\frac{d+1-2dp}{\binom{2d}{d}}\bigg)^{2^d}\prod_{x\in\Z^d\setminus V}2^{-|Y_x|}=\bigg(\frac{d+1-2dp}{\binom{2d}{d}}\bigg)^{2^d}2^{d2^d-2|\partial_ES|}\,,
\end{equation}
which is~\eqref{eq:cut}.
\end{proof}

\paragraph{Counting boundaries.}
To sum~\eqref{eq:cut} over the possible shapes of a finite cluster, we count the possible boundaries. For $m\in\N$ let $\mathcal S_m$ denote the family of finite sets $S\subset\Z^d$ with $o\in S$ and $|\partial_ES|=m$ such that both $S$ and $\Z^d\setminus S$ are connected. Since~\eqref{eq:cut} provides the factor $4^{-m}$, any bound on $|\mathcal S_m|$ with exponential growth rate smaller than $4$ suffices for our purposes. The following elementary bound is enough; a bound with base $3$ goes back to Ruelle~\cite{Ruelle69}, see~\cite[Introduction]{ABE17}, and much better bounds are known in high dimensions~\cite{LM98}.

\begin{lemma}\label{lem:contours}
For every $m\in\N$, we have $|\mathcal S_m|\le m(7/2)^{m-1}/2$.
\end{lemma}
\begin{proof}
Let $S\in\mathcal S_m$ and $F:=\partial_ES$. Every unit square in $\Z^d$, that is, every cycle of length four, contains an even number of edges of $F$, since a walk around the square enters $S$ as often as it leaves $S$. Moreover, $F$ is a minimal edge cut: if $U\subset\Z^d$ satisfies $\emptyset\neq\partial_EU\subset F$, then each of the connected sets $S$ and $\Z^d\setminus S$ lies in $U$ or in $\Z^d\setminus U$, since none of its internal edges belongs to $F$, and they do not lie in the same one, since $\partial_EU\neq\emptyset$; hence $\partial_EU=F$.

We first count the sets $F$ that contain a given edge $f$. Let $F'\subset F$ with $f\in F'$ and $F'\neq F$. Then some unit square contains an odd number of edges of $F'$. Otherwise, every closed walk in $\Z^d$ would traverse edges of $F'$ an even number of times: a closed walk can be reduced to the trivial walk by exchanging consecutive steps in different coordinate directions, which changes the walk along a unit square, and by deleting pairs of consecutive opposite steps, and neither operation changes the parity of the number of steps along edges of $F'$. Hence this parity would be the same for all walks from $o$ to a given vertex $x$, and $F'=\partial_EU$ for the set $U$ of vertices $x$ for which it is even. As $\emptyset\neq F'\subsetneq F$, this contradicts the minimality of $F$.

Fix an enumeration of the unit squares and an order of the edges of every unit square. Starting from $F'=\{f\}$, we repeatedly take the first unit square that contains an odd number of edges of $F'$, and add to $F'$ all edges of this square that belong to $F\setminus F'$; we stop when $F'=F$. As the square contains an even number of edges of $F$, it contains either one edge of $F'$ and one or three further edges of $F$, or three edges of $F'$ and one further edge of $F$. We record a step that adds a single edge by one of three letters of cost one, indicating which of the at most three edges of the square outside $F'$ is added, and a step that adds three edges by a fourth letter of cost three. The resulting word has total cost $|F|-1=m-1$. Since the square in each step is determined by the current set $F'$, the edge $f$ and the word determine $F$. For $\ell\in\N_0$ let $w_\ell$ be the number of words of cost $\ell$. Then $w_0=1$, $w_1=3$, $w_2=9$ and $w_\ell=3w_{\ell-1}+w_{\ell-3}$ for $\ell\ge3$, so by induction $w_\ell\le(7/2)^\ell$, since $3(7/2)^{\ell-1}+(7/2)^{\ell-3}\le(7/2)^\ell$. Hence at most $(7/2)^{m-1}$ of the sets $\partial_ES$, $S\in\mathcal S_m$, contain $f$.

It remains to choose $f$. If $j\ge0$ is minimal with $(j+1)e_1\notin S$, then $\{je_1,(j+1)e_1\}\in F$. Each of the $j+1$ lines $ie_1+\Z e_2$, $i\in\{0,\dots,j\}$, meets the finite set $S$ and therefore contains two edges of $F$; together with the edge $\{je_1,(j+1)e_1\}$, this gives $2(j+1)<m$. Thus $F$ contains one of at most $m/2$ fixed edges. Since $S$ is the connected component of $o$ in the graph obtained from $\Z^d$ by deleting the edges of $F$, the set $S$ is determined by $F$, and the claim follows.
\end{proof}

Combined with the factor $4^{-|\partial_ES|}$ in~\eqref{eq:cut}, Lemma~\ref{lem:contours} yields the summable factor $m(7/8)^m/7$. This is the Peierls argument behind the following estimate.

\begin{proposition}\label{prop:tail}
For all $n\in\N$,
\begin{equation}\label{eq:tail}
\P_p\big(|\calC^\textsc{U}(o)|<\infty,\ \calC^\textsc{U}(o)\not\subset\Lambda_{n-1}\big)\le10n\,(d+1-2dp)^{2^d}\,(7/8)^{2n}\,.
\end{equation}
\end{proposition}
\begin{proof}
Assume that $C:=\calC^\textsc{U}(o)$ is finite and not contained in $\Lambda_{n-1}$, and let $S$ be the union of $C$ and of all finite connected components of $\Z^d\setminus C$. Then $S$ is finite and contains $o$. It is connected, since every finite connected component of $\Z^d\setminus C$ has a neighbor in $C$, and $\Z^d\setminus S$ is the unique infinite connected component of $\Z^d\setminus C$, hence connected. An edge of $\partial_ES$ joins a vertex of $\Z^d\setminus S$ to a vertex of $S$ that does not belong to a finite connected component of $\Z^d\setminus C$, that is, to a vertex of $C$. Hence $\partial_ES\subset\partial_EC$, and $\partial_ES$ is closed.

Moreover, $|\partial_ES|\ge2n+2$. Indeed, $S$ contains $o$ and a vertex $x$ with $|x_i|\ge n$ for some $i\in\{1,\dots,d\}$; fix $j\in\{1,\dots,d\}\setminus\{i\}$. Since $S$ is connected, for every integer $\ell$ between $0$ and $x_i$ there is a vertex $z\in S$ with $z_i=\ell$. The line $z+\Z e_j$ meets the finite set $S$, so it contains at least two edges of $\partial_ES$, and both endpoints of these edges have $i$-th coordinate $\ell$. This yields at least $2(n+1)$ distinct edges of $\partial_ES$.

Thus $S\in\mathcal S_m$ for some $m\ge2n+2$, and $\partial_ES$ is closed. Summing~\eqref{eq:cut} over these sets, and using Lemma~\ref{lem:contours} and $2^d\le\binom{2d}{d}$, we obtain
\begin{equation}
\begin{aligned}
\P_p\big(|\calC^\textsc{U}(o)|<\infty,\ \calC^\textsc{U}(o)\not\subset\Lambda_{n-1}\big)&\le(d+1-2dp)^{2^d}\sum_{m\ge2n+2}\frac m2\Big(\frac72\Big)^{m-1}4^{-m}\\
&=\frac78(2n+9)\,(d+1-2dp)^{2^d}\Big(\frac78\Big)^{2n}\,,
\end{aligned}
\end{equation}
which implies~\eqref{eq:tail}, since $7(2n+9)/8\le10n$.
\end{proof}

\paragraph{Smoothness.}
We now decompose the probability that the origin lies in a finite cluster according to the smallest centered box containing that cluster, and then differentiate this series expansion term by term.

\begin{proposition}\label{prop:factorization}
There is a function $H_d\colon[1/2,(d+1)/(2d)]\to[0,\infty)$ with $H_d((d+1)/(2d))>0$ that is infinitely differentiable, with one-sided derivatives at the endpoints, such that
\begin{equation}\label{eq:factorization}
1-\theta^\textsc{U}_d(p)=(d+1-2dp)^{2^d}H_d(p)\qquad\text{for all }p\in[1/2,(d+1)/(2d)]\,.
\end{equation}
\end{proposition}
\begin{proof}
For $n\in\N$ let $E_n:=\{\calC^\textsc{U}(o)\subset\Lambda_n,\ \calC^\textsc{U}(o)\not\subset\Lambda_{n-1}\}$ be the event that $\Lambda_n$ is the smallest centered box containing $\calC^\textsc{U}(o)$. Since $o$ has at least one out-neighbor, $\calC^\textsc{U}(o)\neq\{o\}=\Lambda_0$. Hence the events $E_n$, $n\in\N$, partition the event $\{|\calC^\textsc{U}(o)|<\infty\}$, and
\begin{equation}\label{eq:order-sum}
1-\theta^\textsc{U}_d(p)=\sum_{n\in\N}\P_p(E_n)\,.
\end{equation}

\emph{Step 1: Polynomial structure of each term.} Fix $n\in\N$. The out-sets $(A_x)_{x\in\Lambda_{n+1}}$ determine all open edges with an endpoint in $\Lambda_n$, hence the open cluster of $o$ in $\Lambda_n$ and whether an open edge leaves it; if none does, this cluster is $\calC^\textsc{U}(o)$. Thus $E_n$ is determined by $(A_x)_{x\in\Lambda_{n+1}}$, and $\P_p(E_n)$ is the sum of the probabilities $\P_p(A_x=B_x\text{ for all }x\in\Lambda_{n+1})$ over the finitely many families $(B_x)_{x\in\Lambda_{n+1}}$ of possible out-sets that realize $E_n$. Each of these probabilities is a product of $|\Lambda_{n+1}|$ factors, each equal to $(d+1-2dp)/\binom{2d}{d}$ or $(2dp-d)/\binom{2d}{d+1}$. If the family $(B_x)$ realizes $E_n$, then the cluster $C:=\calC^\textsc{U}(o)\subset\Lambda_n$ is finite, and every $x\in C$ satisfies $B_x\subset C$ and $|B_x|\ge d$, so it has at least $d$ neighbors in $C$. By Lemma~\ref{lem:corners}, at least $2^d$ vertices $x\in C$ have exactly $d$ neighbors in $C$, and for each of them $B_x$ consists of these $d$ neighbors. Hence each of the products contains at least $2^d$ factors $(d+1-2dp)/\binom{2d}{d}$, and there is a polynomial $h_n$ of degree at most $|\Lambda_{n+1}|$ such that
\begin{equation}\label{eq:order-poly}
\P_p(E_n)=(d+1-2dp)^{2^d}h_n(p)\qquad\text{for all }p\in[1/2,(d+1)/(2d)]\,.
\end{equation}

\emph{Step 2: Termwise differentiation.} Let $J:=[1/2,(d+1)/(2d)]$. Since $E_n\subset\{|\calC^\textsc{U}(o)|<\infty,\ \calC^\textsc{U}(o)\not\subset\Lambda_{n-1}\}$, Proposition~\ref{prop:tail} and~\eqref{eq:order-poly} give $0\le h_n(p)\le10n(7/8)^{2n}$ for $p\in J\setminus\{(d+1)/(2d)\}$, and by continuity for all $p\in J$. By Markov's inequality for derivatives of polynomials~\cite[Theorem~5.1.8]{BorweinErdelyi95}, every real polynomial $P$ of degree at most $M$ satisfies $\max_{[-1,1]}|P'|\le M^2\max_{[-1,1]}|P|$. By an affine change of variables, $\max_J|P'|\le4dM^2\max_J|P|$, since $J$ has length $1/(2d)$. As $h_n$ and its derivatives are polynomials of degree at most $|\Lambda_{n+1}|=(2n+3)^d$, applying this bound $r$ times yields
\begin{equation}
\max_{p\in J}\big|h_n^{(r)}(p)\big|\le\big(4d(2n+3)^{2d}\big)^r\,10n\,(7/8)^{2n}\qquad\text{for all }r\in\N_0\,.
\end{equation}
The right-hand side is summable in $n$, so for every $r\in\N_0$ the series $\sum_{n\in\N}h_n^{(r)}$ converges uniformly on $J$. Consequently, $H_d:=\sum_{n\in\N}h_n$ is nonnegative and infinitely differentiable on $J$, with one-sided derivatives at the endpoints, and~\eqref{eq:factorization} follows from~\eqref{eq:order-sum} and~\eqref{eq:order-poly}.

\emph{Step 3: Positivity.} As in the proof of Corollary~\ref{cor:non-analytic}, we consider the event that $\calC^\textsc{U}(o)$ is the unit cube $\{0,1\}^d$. It occurs if every vertex of the cube selects exactly its $d$ neighbors in the cube and no vertex outside the cube selects a vertex of the cube. Exactly $d 2^d$ vertices outside the cube have a neighbor in the cube, each of them exactly one, and every directed edge is open with probability $p$. By independence,
\begin{equation}
1-\theta^\textsc{U}_d(p)\ge\P_p\big(\calC^\textsc{U}(o)=\{0,1\}^d\big)\ge\bigg((d+1-2dp)\binom{2d}{d}^{-1}\bigg)^{2^d}(1-p)^{d2^d}\,.
\end{equation}
Dividing by $(d+1-2dp)^{2^d}$ and letting $p\uparrow(d+1)/(2d)$, we obtain from~\eqref{eq:factorization} and the continuity of $H_d$ that $H_d((d+1)/(2d))\ge\binom{2d}{d}^{-2^d}\big((d-1)/(2d)\big)^{d2^d}>0$.
\end{proof}

\begin{proof}[Proof of Theorem~\ref{thm:order}]
By Lemma~\ref{lem:forward-infinite} and Proposition~\ref{prop:factorization}, the function $1-\theta^\textsc{U}_d$ vanishes on $[(d+1)/(2d),1)$ and coincides on $(1/2,(d+1)/(2d)]$ with
\begin{equation}
g(p):=(d+1-2dp)^{2^d}H_d(p)=(2d)^{2^d}\big((d+1)/(2d)-p\big)^{2^d}H_d(p)\,,
\end{equation}
which is infinitely differentiable up to $(d+1)/(2d)$. This proves the first claim of~(i), and~(ii) follows with $c:=(2d)^{2^d}H_d((d+1)/(2d))>0$, since $H_d$ is differentiable at $(d+1)/(2d)$. As $2^d$ is even, the Leibniz rule shows that the left derivatives of $g$ at $(d+1)/(2d)$ of orders $1,\dots,2^d-1$ vanish, while that of order $2^d$ equals $(2^d)!\,c$. Since all right derivatives of $1-\theta^\textsc{U}_d$ at $(d+1)/(2d)$ vanish, the two branches join in class $C^{2^d-1}$, which proves~(i) and the first claim of~(iii). The derivatives of order $2^d$ of the two branches of $\theta^\textsc{U}_d$ converge to $-(2^d)!\,c<0$ and to $0$, respectively. In particular, the derivative of order $2^d-1$ of $\theta^\textsc{U}_d$ has distinct one-sided derivatives at $(d+1)/(2d)$, so $\theta^\textsc{U}_d$ is not $2^d$ times differentiable there.
\end{proof}

\subsection*{Disclosure of AI use.}
This project began in January 2026 as a classical research project carried out without AI
tools. Before any AI tools were used, we had convinced ourselves that a
Burton--Keane argument yields uniqueness for $p<(d+1)/(2d)$, and we had
developed the proof of full coverage and local uniqueness for
$p>(d+1)/(2d)$ stated in Theorem~\ref{thm:full-coverage-local-uniqueness}. We also expected
the continuity of the percolation function (Theorem~\ref{thm_con}) to
follow from uniqueness without much additional effort. At a later stage we
used the AI tools Claude Fable and Claude Opus (Anthropic) and ChatGPT Astra
(OpenAI) to work out the details of these proofs, based on our notes. The observation that
a more efficient modification of the Burton--Keane argument also covers the full parameter regime $(1/2d, 1]$ was
suggested by AI when reviewing a preliminary version of this manuscript. The proof of the order of the full-coverage transition (Theorem~\ref{thm:order}) was developed together with AI. We also used these tools for final checks of the manuscript. All authors have checked every argument in the paper and take full responsibility for its content.

\subsection*{Acknowledgements}
BJ and JK gratefully acknowledge funding from the INRIA Team POPOPOP in Lille for a research visit in January 2026 and thank DC and BH for their great hospitality.

\bibliographystyle{alpha}
\bibliography{references}

\end{document}